\documentclass[11pt]{article}
\usepackage[unicode]{hyperref}
\usepackage{amsmath}
\usepackage{graphicx}								
\usepackage{subfig}
\usepackage{amsfonts} 							
\usepackage{dsfont}
\usepackage{bigints} 								
\usepackage{algorithm}              
\usepackage{algorithmic}
\usepackage{booktabs} 							
\usepackage{multirow}								
\usepackage{setspace}
\usepackage{mathrsfs}
\usepackage{upgreek}
\usepackage{tikz}
\usepackage{xfrac}
\usepackage{geometry}

\usepackage[T1]{fontenc}

\newtheorem{theorem}{Theorem}

\newtheorem{Proposition}{Proposition}[section]

\newtheorem{remark}{Remark}

\numberwithin{equation}{section}

\title{Estimates of the distance to the exact solution
of evolutionary reaction-diffusion problems based on local Poincar\'e type inequalities}
\author{S. Matculevich and S. Repin}
\begin{document}
\date{}
\maketitle

\begin{abstract}
The goal of the paper is to derive two-sided bounds of the distance between the exact solution of the
evolutionary reaction-diffusion problem with mixed Dirichlet--Robin boundary conditions
and any function in the admissible energy space. The derivation is based upon
transformation of the integral identity, which defines the generalized solution,
and exploits classical Poincar\'e inequalities and Poincar\'e type inequalities for
functions with zero mean boundary traces. The corresponding constants are estimated
due to \cite{PayneWeinberger1960} and \cite{NazarovRepin2013}. To handle problems with
complex domains and mixed boundary conditions, domain decomposition
is used. The corresponding bounds of the distance to the exact solution,
contain only constants in local Poincar\'e type inequalities associated with subdomains. 
Moreover, it is proved
that the bounds are equivalent to the energy norm of the error.
\end{abstract}

\section{Problem statement}
\label{sc:problem-statement}
We consider the evolutionary reaction-diffusion problem, which is presented in the
following form: find $u = u(x, t)$ and $p = p(x, t)$ such that
\begin{alignat}{3}
	u_t - \nabla \cdot p + \varrho^2 u & =\, f
	&& \quad \mbox{in}\quad Q_T:= \Omega \times (0, T),
  \label{eq:parabolic-problem-equation}\\
  p & =\, A \nabla u
	&& \quad \mbox{in}\quad Q_T, \label{eq:dual-part}\\
	u(\cdot, 0) & =\, u_0
	&& \quad \mbox{in}\quad \Omega,
  \label{eq:parabolic-problem-initial-condition}\\
  u & =\, 0
	&& \quad \mbox{on}\quad S_D,
  \label{eq:parabolic-problem-dirichlet-boundary-condition}\\
  p \cdot n + \sigma^2 u & =\, F
	&& \quad \mbox{on}\quad S_R.
    \label{eq:parabolic-problem-robin-boundary-condition}
\end{alignat}
Here, $\Omega$ is a bounded connected domain in ${\mathds{R}}^d$ ($d \geq 1$) with Lipchitz
continuous boundary $\partial \Omega$, which consists of two measurable non-intersecting
parts $\Gamma_D$ and $\Gamma_R\not=\emptyset$  associated with the Dirichlet and Robin
boundary conditions, respectively, and $n$ denotes the unit outwards normal vector to
$\partial\Omega$. $T $ is a finite positive number,
$S := \partial\Omega \times ] 0, T [$ is the lateral surface of the space--time cylinder
$Q_T$, $S_D: = \Gamma_D \times ] 0, T [$, and $S_R := \Gamma_R \times ] 0, T[$.
%
We assume that $A$ is a symmetric matrix with coefficients in $L^{\infty}(\Omega)$, which  for almost all $x\in \Omega$ satisfies the condition
\begin{equation}
	\underline{\lambda}_A |\xi|^2
	\leq \, A(x) \: \xi \cdot \xi\,
	\leq \, \overline{\lambda}_A | \xi |^2,
	\quad
	\forall \xi \in {\mathds{R}}^d,
	\quad 0 < \underline{\lambda}_A\leq\, \overline{\lambda}_A < +\infty.
	\label{eq:operator-a}
\end{equation}
Here, $|\xi| := \sqrt{\xi \cdot \xi}$, and $\cdot$ stands for the scalar product in ${\mathds{R}}^d$.
Also, we assume that $f\in L^{2}(Q_T)$, $u_0\in L^{2}(\Omega)$, $F \in L^{2}(S_R)$,
and the coefficients $\varrho(x)$ and $\sigma(x)$ are uniformly bounded by the constants
$C_{\varrho}$ and $C_{\sigma}$ in $Q_T$ and $S_T$, respectively.
%
%
%


Throughout the paper, the norms of $L^{2}(\Omega)$ and $L^{2}(Q_T)$ are denoted
by $\| \cdot \|_\Omega$ and $\| \cdot \|_{Q_T}$, respectively.
%
%
$H^{1}_{0}(\Omega)$ is a subspace of $H^1(\Omega)$ with functions
vanishing on Dirichlet part of the boundary, and $V_0 :=$
$H^{1}_{0}(Q_T) := H^1\left(0, T; H^{1}_{0}(\Omega) \right)$.
%
The generalized solution of
(\ref{eq:parabolic-problem-equation})--(\ref{eq:parabolic-problem-robin-boundary-condition})
is a function $u \in L^{2}\left(0, T; H^{1}_{0}(\Omega) \right)$ satisfying the integral identity
\begin{multline}
    \int\limits_{\Omega} \Big( u(x, T) \eta(x, T) - u(x, 0) \eta(x, 0) \Big) {\mathrm{\:d}x}
		- \,\int\limits_{Q_T} u \eta_t {\mathrm{\:d}x\mathrm{d}t}
		+ \,\int\limits_{Q_T} A \nabla{u} \cdot \nabla{\eta} {\mathrm{\:d}x\mathrm{d}t} \\
    + \int\limits_{Q_T} \varrho^2 u \eta {\mathrm{\:d}x\mathrm{d}t}
		+ \,\int\limits_{S_R} \sigma^2 u \eta {\mathrm{\:d}s\mathrm{d}t}
		= \,\int\limits_{Q_T} f \eta {\mathrm{\:d}x\mathrm{d}t} + \int\limits_{S_R} F \eta {\mathrm{\:d}s\mathrm{d}t},
		\quad \forall \eta \in V_0.
    \label{eq:generalized-statement}
\end{multline}
%
Due to known results (see, e.g., \cite{Evans2010, Ladyzhenskaya1985,
Ladyzhenskayaetall1967, Wloka1987}), the solution of (\ref{eq:generalized-statement})
exists and is unique.

Assume that $v \in V_0$ is a function compared with $u$ (e.g., it can be an
approximation generated by some numerical method). Our goal is to deduce explicitly
computable and realistic estimates of the error $e := u - v$ in terms of the measure

%
\begin{equation}
	[ e ]^2_{(\nu, \theta, \zeta, \chi)}
	:= \nu \! \int\limits_0^T \left \| \, \nabla e \, \right\|^2_{A, \, \Omega}  \mathrm{\:d}t \,
	+ \, \int\limits_0^T \left\| \, \theta \, e \right\|^2_\Omega \mathrm{\:d}t \,
	+ \, \zeta \! \, \left \| \, e (\cdot, T) \, \right \|^2_{\Omega} \,
	+ \, \chi \! \int\limits_0^T \left\| \, \sigma\, e \, \right\|^2_{\Gamma_R} \mathrm{\:d}t,
	\label{eq:energy-norm-for-reaction-diff-evolutionary-problem}
\end{equation}
where $\| \nabla e \|^2_{A, \,\Omega} := \int\limits_{\Omega} A  \nabla e \cdot  \nabla e \mathrm{\:d}x$
is the spatial error norm defined on $\Omega$,
%
and $\nu$, $\theta$, $\zeta$, $\chi$ are positive weights. The first two terms present
a measure equivalent to the natural energy norm, the third term measures the error at
$t=T$, and the last one measures possible violations of the Robin boundary condition.
The weights can be selected in order to balance different components of the error in
desired proportions. In other words, the quantity
(\ref{eq:energy-norm-for-reaction-diff-evolutionary-problem}) generates a collection
of different error measures, which can be used for judging on the distance between $u$
and $v$.

In this paper, we derive computable  majorants of $[ e ]^2_{(\nu, \theta, \zeta, \chi)}$
by means of the method close to that has been originally suggested in  \cite{Repin2002}
(see also Section 9.3 of the monograph \cite{RepinDeGruyter2008}). They are derived
by special transformations of the integral identity (\ref{eq:generalized-statement}).
However, in this paper we apply a somewhat new approach based on domain decomposition
and local Poincar\'e type inequalities for functions with zero mean boundary traces
(sharp constants in these inequalities has been recently found in
\cite{NazarovRepin2013}). As a result, we obtain fully computable estimates, which are
applicable for problems with complicated geometry and non-trivial boundary conditions.



In Section \ref{sc:general-estimates}, we deduce two--sided bounds of the distance to the exact
solution. The respective
error majorant (Theorems \ref{th:theorem-minimum-of-majorant-I}),
contains the constants $C_{\mathrm{F\Omega}}$ and $C_{\rm \Gamma R}$ in the Friedrichs and trace
type inequalities
\begin{alignat}{2}
\| & w \|_{\Omega} \leq C_{\mathrm{F\Omega}} \| \, \nabla w \,\|_{\Omega},
\label{eq:friedrichs-inequality}\\
\| & w \|_{\Gamma_R} \leq C_{\rm \Gamma R} \| \,\nabla w \,\|_{ \Omega} \, ,
\label{eq:trace-inequality}
\end{alignat}
which are valid for functions in $H^1_0(\Omega)$.

In general, finding global constants $C_{\mathrm{F\Omega}}$ and $C_{\rm \Gamma R}$ may be not an easy task
(especially for geometrically complicated domains and mixed boundary conditions).
A way to overcome these difficulties is suggested in Section
\ref{sc:domain-decomposition}, where we deduce new forms of majorants, which are based
on decomposition of $\Omega$ into a collection of non-overlapping convex sub-domains
and local constants associated with these subdomains. The constants associated with
subdomains are defined either by means of the  Payne--Weinberger
estimate \cite{PayneWeinberger1960} related to the Poincar\'e inequality for convex
domains or due to the results of \cite{NazarovRepin2013}, where sharp constants for
Poincar\'e type inequalities for functions with zero mean boundary traces has been
found. Therefore, we obtain different forms of the respective error majorants, which
involve only local constants and known functions. In Section 4, we  prove that the majorants are
equivalent to the distance to the exact solution
measured either in terms of the measure
(\ref{eq:energy-norm-for-reaction-diff-evolutionary-problem}) or 
in terms of a combined primal-dual energy norm.

\section{Estimates based on global constants}
\label{sc:general-estimates}

\subsection{Majorant of $[ e ]_{(\nu, \theta, \zeta, \chi)}$}

Let $v \in V_0$ be a function considered as an approximation of $u$. We
transform (\ref{eq:generalized-statement}) and arrive at the relation
\begin{multline}
	\int\limits_0^T \|\,\nabla{e}\,\|^2_{A, \, \Omega} \mathrm{\:d}t \,
	+ \,\int\limits_0^T \left \| \, \varrho \, e  \,\right \|^2_{\Omega} \mathrm{\:d}t \,
	+ \,\tfrac{1}{2} \| \, e(\cdot, T) \,\|^2_{\Omega} \, \\[-10pt]
	+ \, \int\limits_0^T \left \| \, \sigma\, e  \,\right \|^2_{\Gamma_R} \mathrm{\:d}t \,
	= \int\limits_{Q_T} \!\! \left(f - v_t - \varrho^2 \, v \right) \, e {\mathrm{\:d}x\mathrm{d}t} \, \\[-5pt]
	- \,\int\limits_{Q_T} A \nabla{v} \cdot \nabla e {\mathrm{\:d}x\mathrm{d}t} \,
	+ \,\int\limits_{S_R} ( g - \sigma^2 v ) \, e {\mathrm{\:d}s\mathrm{d}t} \, + \,
	\tfrac{1}{2} \| \, e(\cdot  , 0) \,\|^2_{\Omega} \,,
	\label{eq:energy-balance-equation}
\end{multline}
%
which can be viewed as the basic error identity. In order to rearrange
the right hand side of (\ref{eq:energy-balance-equation}), we introduce a vector valued
function $y \in  Y_{{\mathrm{div} \:}}(Q_T)$, where $Y_{\mathrm{div} \:}(Q_T)$ denotes the space of vector
valued functions $y \in L^{2}(\Omega, {\mathds{R}}^d)$ 
$\mathrm{div} \: y \in L^{2}(\Omega)$ and $y \cdot n \in L^{2}(\Gamma_R)$ for almost all
$t \in (0, T)$.  We introduce the quantities
\begin{alignat}{2}
	\mathbf{r}_f  (v, y) & := f - v_t - \varrho^2 \, v + \mathrm{div} \: y, \label{eq:r-f}\\
	\mathbf{r}_A  (v, y) & := y - A \nabla{v}, \label{eq:r-d} \\
	\mathbf{r}_F  (v, y) & := F - \sigma^2 v - y \cdot n, \label{eq:r-b}
\end{alignat}
 which have clear meaning: they are {\em residuals} of
(\ref{eq:parabolic-problem-equation}), (\ref{eq:dual-part}),  and
(\ref{eq:parabolic-problem-robin-boundary-condition}), respectively. The theorem below
shows that certain norms of these quantities control the distance between $u$ and $v$.
Also, we define weighted residuals
\begin{alignat}{2}
	\mathbf{r}_{f,\,\mu} (v, y) & := \mu \, \mathbf{r}_f \quad{\rm and} \quad
	\mathbf{r}_{f, 1 - \mu}  (v, y) := (1 - \mu) \, \mathbf{r}_f.
	\label{eq:r-mu}
\end{alignat}
Here, $\mu(x,t)$ is a real-valued function taking values in $[0, 1]$
(these weighted quantities are motivated later).
\begin{theorem}
\label{th:theorem-minimum-of-majorant-I}
For any $v \in V_0$, $y \in Y_{\mathrm{div} \:}(Q_T)$, $\delta \in (0, 2]$, and
real-valued function $\gamma(t) \geq \frac12$, we have the estimate
\begin{multline}
	{[e]}^{\,2}_{({\nu},\, {\theta},\, 1, \,2)}
	\leq {\overline{\mathrm M}^{\,2}_{\,\mathrm{I}}} (v, y; \delta, \gamma, \mu) \!
	:= \| e (\cdot, 0) \|^2_{\Omega} \,
	+ \, \int\limits_0^T \! \! \Big ( \alpha_1(t) \| \, \mathbf{r}_A (v, y) \,\|^2_{A^{-1}, \Omega} \\[-2pt]
	+ \, \gamma(t) \left\| \, \tfrac{1}{\varrho} \, \mathbf{r}_{f,\, \mu} (v, y) \, \right\|^2_{\Omega} \,
	+ \, \alpha_2(t) \, \tfrac{C_{\mathrm{F\Omega}}^2}{\,\underline{\lambda}_A}
		\| \, \mathbf{r}_{f, 1 - \mu} (v, y) \, \|^2_{\Omega} \,
	+ \, \alpha_3(t) \, \tfrac{C_{\rm \Gamma R}^2}{\,\underline{\lambda}_A}
		\| \, \mathbf{r}_F (v, y) \, \|^2_{\Gamma_R} \Big ) \mathrm{\:d}t,
	\label{eq:majorant-1}
\end{multline}
%
where
%
${\nu} = 2 - \delta$,
${\theta}(x,t) =
\varrho(x) \left( 2 - \tfrac{1}{\gamma(t)} \right)^{\footnotesize \frac{1}{2}}$, and
$\alpha_1(t)$, $\alpha_2(t)$, $\alpha_3(t)$ are arbitrary positive functions satisfying
the relation
\begin{equation}
	\tfrac{1}{\alpha_1(t)} + \tfrac{1}{\alpha_2(t)} + \tfrac{1}{\alpha_3(t)} = \delta.
	\label{eq:alpha}
\end{equation}
\end{theorem}
\noindent
{\bf Proof.}
We transform the right-hand side of (\ref{eq:energy-balance-equation}) by
means of the integral identity
\begin{equation}
    \int\limits_{Q_T} \mathrm{div} \: y \, e {\mathrm{\:d}x\mathrm{d}t} + \int\limits_{Q_T} y \cdot \nabla{e}{\mathrm{\:d}x\mathrm{d}t} =
    \int\limits_{S_R} y \cdot n \,e {\mathrm{\:d}s\mathrm{d}t},
    \label{eq:y-rel}
\end{equation}
which yields
\begin{multline}
    \int\limits_0^T \left \| \, \nabla e \, \right\|^2_{A, \, \Omega}  \mathrm{\:d}t \, + \,
    \int\limits_0^T \left\| \, \varrho \, e \right\|^2_\Omega \mathrm{\:d}t \, + \,
    \tfrac{1}{2} \! \, \left \| \, e (\cdot, T) \, \right \|^2_{\Omega} \, + \,
    \! \int\limits_0^T \left\| \, \sigma\, e \, \right\|^2_{\Gamma_R} \mathrm{\:d}t
    \\ =
    \mathscr{I}_f + \mathscr{I}_A + \mathscr{I}_F + \tfrac{1}{2} \| \, e(x, 0) \, \|^2_{\Omega},
    \label{eq:gen-stat-u-v-norm-y}
\end{multline}
where
\begin{eqnarray}
    &\mathscr{I}_f := &\int\limits_{Q_T} \mathbf{r}_f (v, y) \, e {\mathrm{\:d}x\mathrm{d}t}, \\
    &\mathscr{I}_A := &\int\limits_{Q_T} \mathbf{r}_A (v, y) \cdot \nabla{e} {\mathrm{\:d}x\mathrm{d}t}, \\
    &\mathscr{I}_F := &\int\limits_{S_R} \mathbf{r}_F (v, y) \, e {\mathrm{\:d}s\mathrm{d}t}.
    \label{eq:Ir-Id-Ib-terms}
\end{eqnarray}
%
It is easy to see that
\begin{equation}
    \mathscr{I}_A
    \leq \int\limits_0^T \left \| \, \mathbf{r}_A \, \right\|_{A^{-1}, \, \Omega}
                  \| \nabla{e} \|_{A, \Omega} \mathrm{\:d}t
    \label{eq:id-estimate}
\end{equation}
and (cf. (\ref{eq:trace-inequality}))
\begin{equation}
    \mathscr{I}_F
    \leq \int\limits_0^T \left\| \, \mathbf{r}_F \, \right\|_{\Gamma_R}
                  \|\, e\, \|_{\Gamma_R} \mathrm{\:d}t
		\leq \int\limits_0^T \left\| \, \mathbf{r}_F \, \right\|_{\Gamma_R}
                  \tfrac{C_{\rm \Gamma R}}{\sqrt{\underline{\lambda}_A}}
									\| \,\nabla e\,\|_{A, \Omega} \mathrm{\:d}t.
    \label{eq:ib-estimate}
\end{equation}
%
In order to estimate the term $\mathscr{I}_f$, we apply the same method as in
\cite{RepinSauter2006} and introduce a function $\mu(x, t)$, which takes values in
$[0, 1]$. The idea behind is to split the integral into two parts, which will be later
subject to different parts of the error norm. If the function $\varrho$ has very
different values and may be close to zero, then the resulting estimate is much more
accurate. We can select $\mu$ such that large factors of the type $\tfrac{1}{\varrho}$,
arising in the estimate, are compensated. Hence, we obtain
\begin{equation}
	\mathscr{I}_f
	\leq \int\limits_0^T
	\left ( \left\| \, \tfrac{\mu}{\varrho} \, \mathbf{r}_f \,  \right\|_{\Omega}
					\left \|  \, \varrho \, e  \,  \right \|_{\Omega}
					+ \tfrac{C_{\mathrm{F\Omega}}}{\!\sqrt{\underline{\lambda}_A}} \,
						\left \| \,  (1 - \mu)\, \mathbf{r}_f \,  \right \|_{\Omega}
						\| \nabla{e} \|_{A, \Omega} \right )\mathrm{\:d}t.
  \label{eq:if-estimate}
\end{equation}
%
Combining (\ref{eq:id-estimate})--(\ref{eq:if-estimate}), we find that
\begin{alignat}{2}
	\int\limits_0^T & \left \| \, \nabla e \, \right\|^2_{A, \, \Omega}  \mathrm{\:d}t \,
	+ \, \int\limits_0^T \left\| \, \varrho \, e \right\|^2_\Omega \mathrm{\:d}t \,
	+ \, \tfrac{1}{2} \! \, \left \| \, e (\cdot, T) \, \right \|^2_{\Omega} \,
	+ \, \! \int\limits_0^T \left\| \, \sigma\, e \, \right\|^2_{\Gamma_R} \mathrm{\:d}t \nonumber\\[-10pt]
	& \qquad \leq \tfrac{1}{2} \| \, e(x, 0) \, \|^2_{\Omega}
	+ \int\limits_0^T \Big (
	  \left\| \, \tfrac{\mu}{\varrho} \, \mathbf{r}_f \, \right\|_{\Omega}
		\left \|\, \varrho \, e \right \|_{\Omega}
	  + \tfrac{C_{\mathrm{F\Omega}}}{\!\sqrt{\underline{\lambda}_A}} \,
		  \left \| \, (1 - \mu) \, \mathbf{r}_f \, \right \|_{\Omega}
			\| \nabla{e} \|_{A, \Omega} \nonumber\\[-5pt]
	  & \qquad \qquad \qquad \; \; \:
		+ \left \| \, \mathbf{r}_A \, \right \|_{A^{-1}, \Omega}
		  \| \nabla{e} \|_{A, \Omega}
		+ \tfrac{C_{\rm \Gamma R}}{\!\sqrt{\underline{\lambda}_A}} \,
		  \left\| \, \mathbf{r}_F \, \right\|_{\Gamma_R} \| \nabla{e} \|_{A, \Omega} \Big ) \mathrm{\:d}t.
	\label{eq:estimate}
\end{alignat}
The second term on the right-hand side of (\ref{eq:estimate}) can be estimated by the
Young--Fenchel inequality
%
\begin{align}
	\int\limits_0^T \left\| \,\tfrac{\mu}{\varrho} \, \mathbf{r}_f\,\right\|_{\Omega}
					 \left \| \, \varrho \, e \, \right \|_{\Omega} \mathrm{\:d}t
	\leq
	\int\limits_0^T \tfrac{1}{2} \,
	\left ( \gamma(t)
					\left\| \, \tfrac{\mu}{\varrho} \, \mathbf{r}_f \, \right\|^2_{\Omega}
					+ \tfrac{1}{\gamma(t)} \left \| \, \varrho \, e \, \right \|^2_{\Omega}
	\right) \mathrm{\:d}t,
	\label{eq:young-fenchel-1}
\end{align}
%
%
where $\gamma(t)$ is an arbitrary real-valued function taking values in
$\Big[\tfrac{1}{2}, + \infty\Big[$. Analogously,
\begin{align}
	\int\limits_0^T \tfrac{C_{\mathrm{F\Omega}}}{\!\sqrt{\underline{\lambda}_A}} \,
	         \left \| \, (1 - \mu) \, \mathbf{r}_f \, \right \|_{\Omega}
					 \| \nabla{e} \|_{A, \Omega} \mathrm{\:d}t
  & \leq  \tfrac{1}{2}
	\int\limits_0^T
	\left( \alpha_1(t) \, \tfrac{C_{\mathrm{F\Omega}}^2}{\underline{\lambda}_A}
				 \left \| \, (1 - \mu) \, \mathbf{r}_f  \, \right \|^2_{\Omega}
				+ \tfrac{1}{\! \alpha_1(t)} \| \, \nabla e \, \|^2_{A, \Omega} \right) \mathrm{\:d}t,
	\nonumber
\end{align}

\begin{align}
  \int\limits_0^T \| \, \mathbf{r}_A \, \|_{A^{-1}, \Omega}\, \| \nabla{e} \|_{A, \Omega}  \mathrm{\:d}t &
	\leq \tfrac{1}{2}
	\int\limits_0^T \left( \alpha_2(t) \, \left \| \, \mathbf{r}_A \, \right \|^2_{A^{-1}, \Omega}
	+ \tfrac{1}{ \!\alpha_2(t)} \|\, \nabla e \,\|^2_{A, \Omega} \right) \mathrm{\:d}t,
\end{align}
and
\begin{align}
	\int\limits_0^T
	\left\| \, \mathbf{r}_F \, \right\|_{\Gamma_R}
	\tfrac{ C_{\rm \Gamma R} }{\! \sqrt{\underline{\lambda}_A}} \, \| \nabla{e} \|_{A, \Omega} \mathrm{\:d}t
	\leq
	\tfrac{1}{2}
	\int\limits_0^T \left( \alpha_3(t) \, \tfrac{{C_{\rm \Gamma R}^2}}{\underline{\lambda}_A}
	\left \| \, \mathbf{r}_F \, \right \|^2_{\Gamma_R}
	+ \tfrac{1}{\!\alpha_3(t)} \| \, \nabla e \, \|^2_{A, \,\Omega} \right) \mathrm{\:d}t.
	\label{eq:young-fenchel-4}
\end{align}
%
Here, $\alpha_1(t)$, $\alpha_2(t)$, and $\alpha_3(t)$ are positive functions satisfying
(\ref{eq:alpha}). Then, the estimate (\ref{eq:majorant-1}) follows from
(\ref{eq:young-fenchel-1})--(\ref{eq:young-fenchel-4}).
\hfill $\Box$

\vskip 10pt

\begin{remark}
\rm
The function $y(x,t)$ can be viewed as an approximation of the exact flux
$A\nabla u$.
If it is defined (e.g., by means of some reconstruction of a numerical solution $v$), then the functions
\begin{alignat}{2}
r_1(t) & := \| \, \mathbf{r}_A (v, y) \,\|_{A^{-1}, \Omega}, \\
r_2(t) & := \tfrac{C_{F\Omega}}{\sqrt{\underline{\lambda}_A}}
            \| \, \mathbf{r}_{f, 1 - \mu} (v, y) \, \|_{\Omega}, \\
r_3(t) & := \tfrac{C_{\rm \Gamma R}}{\sqrt{\underline{\lambda}_A}}
            \| \, \mathbf{r}_F (v, y) \, \|_{\Gamma_R}
\end{alignat}
are known. In this case, the majorant ${\overline{\mathrm M}^{\,2}_{\,\mathrm{I}}} (v, y; \delta, \gamma, \mu)$ can be
minimized with respect to the functions $\alpha_1(t)$, $\alpha_2(t)$, $\alpha_3(t)$.
The optimal functions $\alpha^{*}_i(t)$ can be easily found by the method of
Lagrangian multipliers and are defined by the relation
\begin{equation}
\alpha^{*}_i(t) =
\frac{\sum\limits_{i = 1}^{3} r_i(t)}{\, \delta \, r_i(t)}.
\end{equation}
%
However, if we wish to minimize the majorant with respect to $y$, then it is more
advantageous to keep the quadratic structure of (\ref{eq:majorant-1}). In this
case, we can apply iterative minimization procedures similar to those used in
\cite{MaliNeittaanmakiRepin2013, RepinDeGruyter2008} and some other publications cited
therein.
\end{remark}

\begin{remark}
\rm
The majorant $ {\overline{\mathrm M}^{\,2}_{\,\mathrm{I}}} (v, y; \delta, \gamma, \mu)$ has a clear structure. The first term
contains the error in the initial condition and vanishes if the function $v$ exactly
satisfies it. Other terms are formed by norms of the residuals $\mathbf{r}_A$, $\mathbf{r}_f$, and
$\mathbf{r}_F$ and weight factors formed by the global constants $C_{\mathrm{F\Omega}}$ and $C_{\rm \Gamma R}$
related to $\Omega$. Since $v$ satisfies the boundary condition on $\Gamma_D$, the
majorant vanishes if and only if all the residuals are equal to zero, i.e., if and
only if $v$ coincides with $u$ and $y$ coincides with $A\nabla u$.
\end{remark}
\vskip 10pt

\subsection{Minorant of $[ e ]_{(\nu, \theta, \zeta, \chi)}$}

Computable minorants of the deviations from the exact solution of partial differential
equations provide useful information, which allows us to judge on the quality of
error majorants. For elliptic variational problems, a minorant can be derived fairly
easily by means of variational arguments (see \cite{NeittaanmakiRepin2004}). In
\cite{RepinDeGruyter2008}, another derivation method, which does not exploit variational arguments, was suggested. Below, we apply this  method to the considered class of
parabolic problems and deduce computable minorants of the distance to the exact
solution.


\begin{theorem}
\label{th:theorem-mininum-of-minorant}
Let $v, \: \eta \in V_0$, then the following estimate holds:
\begin{equation}
    {\underline{\mathrm M}^2} (v) :=
    \sup\limits_{\eta \, \in \,  V_0} \Bigg \{
    \sum\limits_{i = 1}^{5} G_{i}(\eta,v,\kappa_i) + G_0(\eta,f, F, u_0) \Bigg \} \leq
        {[e]}^{\,2}_{(\underline{\nu},\, \underline{\theta},\, \underline{\zeta}, \, \underline{\chi})}
    \label{eq:lower-estimate}
\end{equation}
where
\begin{alignat}{2}
G_{1}(v,\eta,\kappa_1) & := \int\limits_{Q_T}
\Big( - \nabla \eta \cdot A \nabla v - \tfrac{1}{2 \kappa_1} A\nabla\eta\cdot\nabla\eta  \Big)
{\mathrm{\:d}x\mathrm{d}t}, \quad \nonumber \\[-3pt]
G_{2}(v,\eta,\kappa_2) & := \int\limits_{Q_T}
\Big( \eta_t v - \tfrac{1}{2 \kappa_2}|\eta_t|^2 \Big) {\mathrm{\:d}x\mathrm{d}t}, \quad \nonumber \\[-3pt]
G_{3}(v,\eta,\kappa_3) & := \int\limits_{Q_T}
\varrho^2 \Big( -  v \eta - \tfrac{1}{2 \kappa_3}|\eta|^2 \Big) {\mathrm{\:d}x\mathrm{d}t}, \quad \quad \nonumber \\[-3pt]
G_{4}(v,\eta,\kappa_4) & := \int\limits_{\Omega} \Big( - v(x, T) \eta(x, T) - \tfrac{1}{2 \kappa_4} \, |\eta(x, T)|^2 \Big) \mathrm{\:d}x, \quad \nonumber \\[-3pt]
G_{5}(v,\eta,\kappa_5) & := \int\limits_{S_R} \sigma^2 \Big( - v \eta - \tfrac{1}{2 \kappa_5 }|\eta|^2 \Big) {\mathrm{\:d}s\mathrm{d}t}, \quad\\
G_0(\eta, f, F, u_0)& :=  \int\limits_{Q_T} f \eta {\mathrm{\:d}x\mathrm{d}t} +
                    \int\limits_{S_R} F \eta {\mathrm{\:d}s\mathrm{d}t} + \int\limits_{\Omega} u_0\eta(\cdot, 0) \mathrm{\:d}x,
\end{alignat}
and
%
$\underline{\nu} = \frac{\kappa_1}{\,2},\;
\underline{\theta}(x) = \Big( \frac{1}{2}
\big( \kappa_2 + \kappa_3 \, \varrho(x)^2 \big) \Big)^{\footnotesize \frac{1}{\,2}},\;
\underline{\zeta} = \frac{\kappa_4}{\,2},\;
\underline{\chi} = \frac{\kappa_5}{\,2},$
%
and $\kappa_i$ ($i = 1, \ldots, 5$) are arbitrary positive numbers.
\end{theorem}

\noindent
{\bf Proof.}
Consider the functional
\begin{multline*}
    \!\mathcal{M} (e) := \\
		\sup\limits_{\eta \in V_0}
    \Bigg \{\!
    \int\limits_{Q_T} \! \! \bigg(
    \nabla \eta \cdot A \nabla e - \tfrac{1}{2 \kappa_1}\, A\nabla\eta\cdot\nabla\eta -
        \eta_t e - \tfrac{1}{2 \kappa_2}\, |\eta_t|^2 +
    \varrho^2 \Big( e \eta - \tfrac{1}{2 \kappa_3} \,|\eta|^2 \Big)
    \bigg) {\mathrm{\:d}x\mathrm{d}t} \\[-5pt]
    \qquad \qquad \qquad \qquad
    + \int\limits_{\Omega} \Big( e(x, T) \eta (x, T) - \tfrac{1}{2 \kappa_4} \, |\eta (x, T)|^2 \Big)\mathrm{\:d}x +
    \int\limits_{S_R} \sigma \Big( e \eta  - \tfrac{1}{2 \kappa_5} \, |\eta|^2 \Big){\mathrm{\:d}s\mathrm{d}t} \Bigg \}. 
%
\end{multline*}
It is not difficult to see that for any $\eta \in V_0$
\begin{alignat}{3}
\displaystyle \int\limits_{Q_T} \bigg( \nabla \eta \cdot A \nabla e
- \tfrac{1}{2 \kappa_1} A\nabla\eta\cdot\nabla\eta \bigg) {\mathrm{\:d}x\mathrm{d}t}
& \leq \frac{\kappa_1}{2} \int\limits_0^T  \| \nabla e \|^2_{A} \mathrm{\:d}t,\\
\displaystyle \int\limits_{Q_T}
\Big( - \eta_t e - \tfrac{1}{2 \kappa_2} |\eta_t|^2 \Big) {\mathrm{\:d}x\mathrm{d}t}
& \leq \frac{\kappa_2}{2}\displaystyle \int\limits_0^T \| e \|^2_{\Omega} \mathrm{\:d}t, \\
\displaystyle \int\limits_{Q_T}
\varrho^2 \Big( e \eta - \tfrac{1}{2 \kappa_3} |\eta|^2 \Big) {\mathrm{\:d}x\mathrm{d}t}
& \leq \frac{\kappa_3}{2} \int\limits_0^T \| \varrho \, e \|^2_{\Omega} \mathrm{\:d}t, \\
\displaystyle \int\limits_{\Omega}
\Big( e(x, T) \eta (x, T) - \tfrac{1}{2 \kappa_4} |\eta (x, T)|^2 \Big)\mathrm{\:d}x
& \leq \frac{\kappa_4}{2} \, \| e(x, T) \|^2_{\Omega},  \\
\displaystyle \int\limits_{S_R}
\sigma \Big( e \eta  - \tfrac{1}{2 \kappa_5} |\eta|^2 \Big) {\mathrm{\:d}s\mathrm{d}t}
& \leq \frac{\kappa_5}{2} \int\limits_0^T \| \sigma^{\tfrac12} e \|^2_{\Gamma_R} \mathrm{\:d}t.
\end{alignat}
Hence, we find that
\begin{equation}
    \mathcal{M} (e)
    \leq
        {[e]}^{\,2}_{(\underline{\nu},\, \underline{\theta},\, \underline{\zeta}, \, \underline{\chi})}.
\end{equation}
%
%
By means of (\ref{eq:generalized-statement}) we rewrite $\mathcal{M} (e, \eta)$ in the
form
\begin{equation*}
    \mathcal{M} (e)
		= \sup\limits_{\eta \, \in \,  V_0} \Bigg \{
      \sum\limits_{i = 1}^{5} G_{i}(\eta,v,\kappa_i) + G_0(\eta,f, F, u_0) \Bigg \},
\end{equation*}
where $\eta$ is any function in $V_0$ and, therefore, we arrive at
(\ref{eq:lower-estimate}).
\begin{remark}
\rm
${\underline{\mathrm M}^2} (v)$ vanishes if and only if $v$ coincides with $u$.
\end{remark}

%

\section{Estimates based on domain decomposition and local constants}
\label{sc:domain-decomposition}

\subsection{Estimates of constants in local Poincar\'{e} type inequalities}
\label{sc:constants}

The majorant defined in Theorem \ref{th:theorem-minimum-of-majorant-I} contains global
constants $C_{\mathrm{F\Omega}}$ and $C_{\rm \Gamma R}$. In general, finding these constants may be not an
easy task (which is equivalent to deriving a guaranteed lower bound of the minimal
eigenvalue for the respective differential operator). Below, we suggest the method,
which allows us to overcome this difficulty. The key idea is to decompose $\Omega$
(which may have a complicated structure) into a collection of simple sub-domains and
derive such an estimate of the distance to the exact solution that uses only local
constants associated with subdomains (a close method for elliptic problems is
considered in [15]). 
We note that for the minorant $\underline {\rm M}$ such a procedure is not required because 
it  does not contain  constants related to the Friedrichs or
Poincar\'e inequalities.

Assume that
%
\begin{equation}
    \overline{\Omega} := \bigcup\limits_{ \Omega_i \subset \, \mathcal{O}_\Omega}
    {\overline{\Omega}}_i,
    \quad
    {\Omega}_{i} \, \cap \, {\Omega}_{j} = \emptyset,\;\;\;
                                                    i \neq j, \;
                                                    i,j = 1, \ldots, N,
    \label{eq:omega-representation}
\end{equation}
where $\Omega_i$ are convex domains with Lipschitz boundaries, and $\mathcal{O}_\Omega$ 
is the partition formed by subdomains $\Omega_i$ (in practice $\{\Omega_i\}_{i = 1}^{N}$ are
typically simplicial or polyhedral cells).  Henceforth, we use the notation
$\Gamma_{ij} = \overline{\Omega}_i \cap {\overline{\Omega}_j}$,
$\Gamma_{Di}=\overline{\Omega}_i \cap \Gamma_D$,
and $\Gamma_{Ri} = \overline{\Omega}_i \cap \Gamma_R$.

For any $\Omega_i$ we have the classical Poincar\'e inequality \cite{Poincare1890}
\begin{equation}
    \|w\|_{\Omega_i} \leq C_{P\Omega_i} \|\nabla w\|_{\Omega_i},
    \label{eq:poincare-equation-general}
\end{equation}
which holds for any function
\begin{equation*}
	w \in \widetilde{H}^1 (\Omega_i) :=
	\left\{ \, w \in H^1(\Omega_i)\,  \big | \, { \big \{  w \big\} }_{\Omega_i} = 0 \, \right\},
\end{equation*}
where $\{w\}_{\Omega_i} := \tfrac{1}{|\Omega_i|} \int\limits_{\Omega_i} w \mathrm{\:d}x$.
%
Due to \cite{PayneWeinberger1960}, we know that
$C_{P\Omega_i} \leq \tfrac{{\rm diam}\, \Omega_i}{\pi}$. This estimate of the Poincar\'e constant admits various generalizations (see, e.g., \cite{AcostaDuran2009,
ChuaWheeden2006};
similar estimates for spaces of vector-valued functions are considered in
\cite{Fuchs2011}).

Poincar\'e type estimates also hold for functions having zero mean traces on the
boundary. Let
\begin{equation}
    \widetilde{H}^1(\Omega_i, {\mathcal T}) :=
    \left\{ w \in H^1(\Omega_i)\,  \big | \,{ \big \{  w \big\} }_{\mathcal{T}} =0 \right\},
    \label{eq:space-with-boundary-mean}
\end{equation}
where ${\mathcal T}$ is a part of the boundary $\partial\Omega_i$, which coincides
with $\Gamma_{ij}$ or $\Gamma_{Ri}$.
For any \linebreak $w \in \widetilde{H}^1 (\Omega_i, {\mathcal T})$, we have the estimate
\begin{alignat}{2}
    \|w\|_{{\mathcal T}} & \leq C_{{\mathcal T}\Omega_i} \|\nabla w\|_{\Omega_i}.
    \label{eq:c-omega-gamma-trace-inequality}
\end{alignat}
Sharp values of $C_{{\mathcal T}\Omega_i}$ are found in \cite{NazarovRepin2013} for
some classes of domains. For our subsequent analysis, we need results related to the
cases, where $\Omega_i$ is either a triangle or a quadrilateral and ${\mathcal T}$ is
one side of it.  We can extend these results to the case of $d=3$ and $\Omega_i$
presented by parallelepiped (or domains obtained by affine transformations of
parallelepiped). Below, for the convenience of the reader, we recall some of these
results.
\begin{itemize}
\item[1.]
If $d = 2$, $\Omega_i$ is the right quadrilateral $\Pi_2 := (0, h_1) \times (0, h_2)$,
and ${\mathcal T}$ is the face $x_1 = 0$, then
\begin{equation}
C_{{\mathcal T}\Pi_2} = \bigg( \tfrac{\pi}{h_{2}}
\tanh \left(\tfrac{\pi h_1}{ h_{2}} \right) \bigg)^{\footnotesize \scalebox{0.5}[1.0]{\( - \)} \frac{1}{2}}.
\end{equation}
Analogously, if $d = 3$, $\Pi_3 := (0, h_1) \times (0, h_2)\times(0,h_3)$, and
${\mathcal T}$ is again the face defied by the condition $x_1 = 0$, then
\begin{equation}
	C_{{\mathcal T}\Pi_3} =
	\bigg( \tfrac{\pi}{h_{+}} \tanh \left(\tfrac{\pi h_1}{ h_{+}} \right)
	\bigg)^{\footnotesize \scalebox{0.5}[1.0]{\( - \)} \frac{1}{2}}, \quad
	h_{+} = \max \left\{ h_2, h_3 \right\}.
	\label{eq:rect-angle-d}
\end{equation}

\item[2.]
If $d = 2$, $\Omega_i$ is the triangle
$\overline{T} := {\rm conv} \Big\{ (0, 0), (0, h), (h, 0) \Big\}$, and ${\mathcal T}$
is the leg defined by the condition $x_1 = 0$, then
$ C_{{\mathcal T}T} = \left(\tfrac{h}{\sigma_1} \right)^{\footnotesize \frac12}$,
where
$\sigma_1 = \zeta_1 \, \tanh(\zeta_1)$, and $\zeta_1$ is the unique root of the equation
$\tan(z) + \tanh(z)  = 0$ in $(0, \pi)$.


\item[3.]
Also, we may use another result of \cite{NazarovRepin2013} related to the case, where
functions have zero mean values on the hypotenuse of the isosceles right triangle $T$
with legs $h$. In this case,
$C_{{\mathcal T}T} = \left(\tfrac{h}{2}\right)^{\footnotesize \frac12}$.
\end{itemize}
%
By means of 2 and 3 and standard affine transformation of the coordinates, we can obtain estimates of $C_{{\mathcal T}T}$ for any non-degenerate triangle.
%
%
\begin{Proposition}
\label{pr:on-boundary}
Let $ T $ be the triangle with the nodes
$\Big\{ (0, 0), (h_1,0), (h_2 \cos \alpha, \, h_2\sin \alpha) \Big\}$ and
${\mathcal T} := \Big\{ x_1 \in [0, h_1] ;\;\; x_2 = 0 \Big\}$. Then, for any
$v\in H^1(T)$ with zero mean trace on $\mathcal T$ we have the estimate
\begin{equation}
	\|v\|_{\mathcal T} \, \leq\,
	C_{{\mathcal T}T} \, h_1^{\footnotesize \frac{1}{2}} \,\|\nabla v\|_T, \quad
	C_{{\mathcal T}T} = \widehat C_{{\mathcal T}T} \, \widehat{C}(\rho, \alpha),
	\label{eq:lemma-triangle}
\end{equation}
where
\begin{equation}
\widehat{C} (\rho, \alpha) =
\Big( \tfrac{\mu(\rho)}{\rho \sin \alpha}\Big)^{\footnotesize \frac{1}{2}}, \, \quad
\mu(\rho) = \tfrac{1}{2}
\Big( 1 + \rho^2 +
\big( 1 + \rho^4 + 2 \cos (2 \alpha) \, \rho^2 \big)^{\footnotesize \frac{1}{2}} \Big), \quad
\rho \,=\, \tfrac{h_2}{h_1},
\end{equation}
and $\widehat C_{{\mathcal T}T} $ is the corresponding constant for the basic
right triangle.
\end{Proposition}
%
\begin{remark}
\rm
It is clear that for the inequality (\ref{eq:c-omega-gamma-trace-inequality}), we have
a certain monotonicity property, which allows us to easily estimate the constant
$C_{\mathrm{\Gamma\Omega}}$.
Namely, if $\Omega_1$ and $\Omega_2$ have a common part $\Gamma$ and
$\Omega_1 \subset \Omega_2$, then
\begin{equation}
	\| w \|
	\leq C_{\mathrm{\Omega_1}} \| \nabla w \|_{\Omega_1}
	\leq C_{\mathrm{\Omega_1}} \| \nabla w \|_{\Omega_2}, \Rightarrow 
	C_{\mathrm{\Gamma\,\Omega_2}} \, \leq \, C_{\mathrm{\Gamma\,\Omega_1}}.
\end{equation}
Therefore,
$C_{\mathrm{\Gamma\,\Omega_2}} \, \leq \, C_{\mathrm{\Gamma\,\Omega_1}}$.
\end{remark}

\subsection{The first estimate}
Let the sub-domains be collected into two different sets
%
\begin{alignat}{2}
	\overline{\Omega}_{\mathrm{P}}
	& := \bigcup\limits_{\Omega_l \subset \, \mathcal{O}_{\mathrm{P}}} {\overline{\Omega}}_l,
	\quad
	\mathcal{O}_{\mathrm{P}} :=
	\Big\{ \; \Omega_l \subset \, \mathcal{O}_\Omega \; \big|
				 \; \varrho|_{\Omega_l} \geq \mathrm{P} , \;
						l = 1, \ldots, N_{\mathrm{P}}  \; \Big\}, \quad
	\mbox{and}
	\label{eq:set-varrho-more-c}
	\\
	\overline{\Omega}_{0}
	& := \bigcup\limits_{\Omega_k \subset \, \mathcal{O}_0} {\overline{\Omega}}_k,
	\quad
	\mathcal{O}_{0} :=
	\Big \{ \; \Omega_k \; \subset \, \mathcal{O}_\Omega \big|
					\; \varrho|_{\Omega_k} < \mathrm{P},\; k = 1, \ldots, N_0  \; \Big \},
	\label{eq:set-varrho-geq-or-less-c}
\end{alignat}
which contain regions with relatively large and small reaction, respectfully.
For the sub-domains in $\mathcal{O}_{0}$, we impose an additional condition, namely,
\begin{equation}
    {{ \Big \{ \mathbf{r}_{f, 1 - \mu} (v, y) \Big\} }_{{\Omega}_k \subset\, \mathcal{O}_0}} = 0, \quad
    \mbox{for a.a.} \quad t \in [0, T].
    \label{eq:mean-condition-on-omega-0}
\end{equation}
%
Since $y$ is in our disposal, then selecting it in such a way that the mean value
condition (\ref{eq:mean-condition-on-omega-0}) holds is technically not difficult.

We impose similar local type conditions on $\Gamma_R$, which is decomposed
into \linebreak
${\Gamma_R}_j = \partial \Omega_j \cap \Gamma_R$, $j = 1, \ldots, M$, $M \leq N$.
%
%
%
%
%
%
%
%
Assume that
%
\begin{equation}
		{{ \Big \{ \mathbf{r}_F (v, y) \Big\} }_{{\Gamma_R}_{j} \subset \mathcal{\mathcal{S}}_{R}}} = 0, \; \;
    \mbox{for a.a.} \quad t \in [0, T],
    \label{eq:mean-condition-on-boundary}
\end{equation}
holds. Here, $\mathcal{\mathcal{S}}_{R}$ denotes a
collection of non-overlapping faces ${\Gamma_R}_{j}$.
%
%
%
%
%
%
Using the idea of Proposition \ref{pr:on-boundary}, we deduce another form of the error majorant, which involves constants in Poincar\'e type inequalities.
Henceforth, we use the following quantities based on localized residuals
\begin{alignat}{2}
    & R_{\mathcal{O}_{\mathrm{P}}, \{ \mathbf{r}_{f, 1 - \mu} \}} (t) :=
    \sum\limits_{ \Omega_l \subset \, \mathcal{O}_{\mathrm{P}} }
    \tfrac{|\Omega_l|}{\!\mathrm{P}^2} \,
    {{ \Big \{  \mathbf{r}_{f, 1 - \mu} (v, y)\Big\} }^2_{\Omega_l}},
		\label{eq:r-op-mean}\\
    & R_{\,\mathcal{O}_{\mathrm{P}}, \, \| \mathbf{r}_{f, 1 - \mu} \|} (t) :=
    \sum\limits_{ \Omega_l \subset \, \mathcal{O}_{\mathrm{P}} }
    \tfrac{{C_{\mathrm{P\Omega}}}_l^2}{\!\underline{\lambda}_A} \, \left\| \mathbf{r}_{f, 1 - \mu} (v, y) \right\|^2_{\Omega_l},
		\label{eq:r-op-norm}\\
    & R_{\, \mathcal{O}_0} (t) \,:=\,
    \sum\limits_{\Omega_k \subset \, \mathcal{O}_0} \tfrac{{C_{\mathrm{P\Omega}}}_k^2}{\!\underline{\lambda}_A} \, \left\|    \mathbf{r}_{f, 1 - \mu} (v, y) \, \right\|^2_{\Omega_k}, \qquad
		\label{eq:r-oo}\\
    & R_{\mathcal{S}_{R}}(t) :=
    \sum\limits_{{\Gamma_R}_{j} \subset \, \mathcal{S}_{R}}
    \tfrac{{C_{\mathrm{\Gamma\Omega}}}_{j}^2}{\!\underline{\lambda}_A} \,\left\| \mathbf{r}_F (v, y) \right\|^2_{{\Gamma_R}_{j}}.
		\label{eq:r-sr}
\end{alignat}
%

\begin{theorem}
\label{th:theorem-majorant-for-decomposed-domain-1}
\noindent
(i) Assume that (\ref{eq:mean-condition-on-omega-0}) and (\ref{eq:mean-condition-on-boundary})
hold,
%
%
then for any $v \in V_0$ and \linebreak
$y \in Y_{\mathrm{div} \:}(Q_T)$, $\delta \in (0, 2]$, $\rho_1(t) \geq 1$, 
$\rho_2(t) \geq 1$, we have the estimate
\begin{multline}
    {[e]}^{\,2}_{({\nu},\, {\theta},\, 1,\, 2)} \leq
    {\overline{\mathrm M}^{\,2}_{\,\rm I, N}} (v, y; \delta, \rho_1, \rho_2, \mu)\! :=
    \int\limits_0^T \! \! \Big (
    \rho_1 \! \left \| \tfrac{1}{\varrho} \, \mathbf{r}_{f,\, \mu} (v, y)  \, \right \|^2_{\Omega} \, + \,
    \rho_2 R_{\mathcal{O}_{\mathrm{P}}, \{ \cdot \}} (t) \, \\[-5pt]
    +
    \alpha_1(t) \| \mathbf{r}_A (v, y) \, \|^2_{A^{-1}, \Omega} \, + \,
    \alpha_2(t) \Big( R_{\mathcal{O}_{\mathrm{P}}, \| \cdot \|}(t) \, + \,
    R_{\, \mathcal{O}_0} (t)    \Big) \, + \,
    \alpha_3(t) R_{\,\mathcal{S}_{R}} (t) \!
    \Big ) \! \mathrm{\:d}t,
\label{eq:majorant-decomposed-i}
\end{multline}
%
where
$\mathbf{r}_f (v, y)$, $\mathbf{r}_{f,\, 1 - \mu}(v, y)$ and $\mathbf{r}_{f,\, \mu}(v, y)$, $\mathbf{r}_A(v, y)$, $\mathbf{r}_F(v, y)$ are defined in (\ref{eq:r-f}), (\ref{eq:r-mu}), (\ref{eq:r-d}), and (\ref{eq:r-b}), respectively,
${\nu} = 2 - \delta$,
${\theta}(x) =
\varrho(x) \left(2 - \tfrac{1}{\rho_1(t)} - \tfrac{1}{\rho_2(t)} \right)^{\footnotesize \frac{1}{\,2}}$
are positive weights,
$\mu(x,t)$ is a real-valued function taking values in $[0, 1]$,
the reaction function $\varrho(x) > 0$,
$\alpha_1(t)$, $\alpha_2(t)$, $\alpha_3(t)$ are positive scalar-valued
functions satisfying the relation (\ref{eq:alpha}).
%
\vskip 10pt
\noindent
(ii) For any $\delta \in (0, 2]$,
$\rho_1(t) \geq 1$, $\rho_2(t) \geq 1$,
and a real-valued function $\mu(x, t)$ taking values in $[0, 1]$,
the upper bound of the variation problem generated by the majorant
\begin{equation}
    \inf\limits_{
    \begin{array}{c}
    v \in V_0\\
    y \in Y_{\mathrm{div} \:}(Q_T)
    \end{array}
    } \overline{\mathrm M}^{\,2}_{\,\rm I, N} (v, y; \delta, \rho_1, \rho_2, \mu)
    \label{eq:inf-maj-I}
\end{equation}
is zero, and it is attained if and only if $v = u$ and $y = A \nabla u$.

\end{theorem}

\noindent
{\bf Proof.}
We consider (\ref{eq:gen-stat-u-v-norm-y}) and estimate $\mathscr{I}_A$ and $\mathscr{I}_F$ analogously
to the proof of \linebreak
Theorem \ref{th:theorem-minimum-of-majorant-I}. The term $\mathscr{I}_f$ is decomposed as
follows:
\begin{multline}
    \mathscr{I}_f \,=\,
    \int\limits_0^T \bigg( \;
    \int\limits_{\Omega} \mathbf{r}_{f,\, \mu} \, e \mathrm{\:d}x \, + \, \int\limits_{\Omega} \mathbf{r}_{f,\, 1 - \mu} \, e \mathrm{\:d}x
    \bigg) \! \mathrm{\:d}t \, \\
    =
    \int\limits_0^T \bigg( \;
    \int\limits_{\Omega} \mathbf{r}_{f,\, \mu} \, e \mathrm{\:d}x \, + \,
    \int\limits_{\Omega_{\mathrm{P}}} \mathbf{r}_{f,\, 1 - \mu} \, e \mathrm{\:d}x  \, + \,
    \int\limits_{\Omega_0} \mathbf{r}_{f,\, 1 - \mu} \, e \mathrm{\:d}x
    \bigg) \! \mathrm{\:d}t \,
    \\
    =
    \mathscr{I}_{f, \,\mu} \, + \, {\mathscr{I}}^{\mathcal{O}_{\mathrm{P}}}_{f, \,1 - \mu} \, + \, {\mathscr{I}}^{\mathcal{O}_0}_{f, \,1 - \mu}.
    \label{eq:representation-of-if-with-mu-maj-1}
\end{multline}
%
Each term on the right-hand side of (\ref{eq:representation-of-if-with-mu-maj-1}) is
estimated by different methods. We use the H\"older inequality, to estimate
$\mathscr{I}_{f, \,\mu}$.
%
%
If (\ref{eq:mean-condition-on-omega-0}) holds, the term
${\mathscr{I}}^{\mathcal{O}_0}_{f,\, 1 - {\mu}}$ can be estimated by
(\ref{eq:poincare-equation-general})
\begin{equation}
    {\mathscr{I}}^{\mathcal{O}_0}_{f,\, 1 - {\mu}} \leq
    \int\limits_0^T R^{\tfrac{1}{2}}_{\, \mathcal{O}_0} \; \| \nabla e \|_{A, \Omega_0} \mathrm{\:d}t.
    \label{eq:omega-0}
\end{equation}
After the following representation
\begin{alignat}{2}
    {\mathscr{I}}^{\mathcal{O}_{\mathrm{P}}}_{f, 1 - \mu} \, & = \,
    \int\limits_0^T \Bigg( \sum\limits_{\Omega_l \subset \,\mathcal{O}_{\mathrm{P}}} \,
		                \int\limits_{\Omega_l} \widetilde{\mathbf{r}}_{f, 1 -    \mu} \, e \mathrm{\:d}x \, + \,
                    \sum\limits_{\Omega_l \subset \,\mathcal{O}_{\mathrm{P}}}
										{ \Big \{ \mathbf{r}_{f, 1 - \mu} \Big\} }_{\Omega_l} \int\limits_{\Omega_l} e \mathrm{\:d}x \Bigg) \mathrm{\:d}t,
\label{eq:majorant-on-decomposed-omega}
\end{alignat}
%
%
the term ${\mathscr{I}}^{\mathcal{O}_{\mathrm{P}}}_{f, \,1 - \mu}$ is estimated as follows
\begin{alignat}{2}
	\int\limits_0^T \sum\limits_{\Omega_l \subset \,\mathcal{O}_{\mathrm{P}}} \;
					 \int\limits_{\Omega_l} \widetilde{\mathbf{r}}_{f, 1 - \mu} e \mathrm{\:d}x \mathrm{\:d}t \, &
	\leq
	\int\limits_0^T R^{\footnotesize \frac{1}{2}}_{\mathcal{O}_{\mathrm{P}}, \| \cdot \|} \;
					 \| \nabla e \|_{A, \Omega_{\mathrm{P}}} \mathrm{\:d}t \,,
	\label{eq:omega-p}
	\\
	%
	\int\limits_0^T \sum\limits_{\Omega_l \subset \,\mathcal{O}_{\mathrm{P}}}
	                { \Big \{ \mathbf{r}_{f, 1 - \mu} \Big\} }_{\Omega_l} \int\limits_{\Omega_l} e \mathrm{\:d}x \mathrm{\:d}t \,  &
	\leq
	\int\limits_0^T \sum\limits_{\Omega_l \subset \,\mathcal{O}_{\mathrm{P}}}
	         \tfrac{|\Omega_l|^{\footnotesize \frac{1}{\,2}}}{\,\mathrm{P}} \,
					 { \Big \{ \mathbf{r}_{f, 1 - \mu} \Big\} }_{\Omega_l} \big \|\varrho \, e \big \|_{\Omega_l} \mathrm{\:d}t \, 		
	\nonumber\\
	& \leq
	\int\limits_0^T R^{\footnotesize \frac{1}{\,2}}_{\mathcal{O}_{\mathrm{P}}, \{ \cdot \}} \;
	\big \|\varrho \, e \big \|_{\Omega} \mathrm{\:d}t.
	\label{eq:omega-p-mean}
\end{alignat}
%
By means of the Minkowski inequality, the sum of right-hand sides of
(\ref{eq:omega-0}) and (\ref{eq:omega-p}) is estimated as follows:
\begin{equation}
\int\limits_0^T R^{\footnotesize \frac{1}{2}}_{\, \mathcal{O}_0} \;
         \| \nabla e \|_{A, \Omega_0} \mathrm{\:d}t +
\int\limits_0^T R^{\footnotesize \frac{1}{2}}_{\mathcal{O}_{\mathrm{P}}, \| \cdot \|} \;
         \| \nabla e \|_{A, \Omega_{\mathrm{P}}} \mathrm{\:d}t \leq
\int\limits_0^T \Big( R_{\, \mathcal{O}_0} +
							R_{\mathcal{O}_{\mathrm{P}}, \| \cdot \|}
				 \Big)^{ \footnotesize \frac{1}{\,2}}
				 \| \nabla e \|_{A, \Omega} \mathrm{\:d}t.
\end{equation}
%
We recall (\ref{eq:mean-condition-on-boundary}) and apply
(\ref{eq:c-omega-gamma-trace-inequality}) to obtain
\begin{equation}
    \mathscr{I}_F \!
    =
    \int\limits_0^T \!\!
   \sum\limits_{{\Gamma_R}_{j} \subset \, \mathcal{S}_{R}} \,
  \int\limits_{{\Gamma_R}_{j}} \widetilde{\mathbf{r}}_F  \, e {\mathrm{\:d}s} \mathrm{\:d}t \leq
    \int\limits_0^T \! R_{\,\mathcal{S}_{R}}^{\footnotesize \frac{1}{\,2}} \| \nabla e \|_{A, \Omega} \mathrm{\:d}t \,.
  \label{eq:estimation-i-b}
\end{equation}
%
In view of the Young--Fenchel inequalities, we have
\begin{alignat}{2}
    \int\limits_0^T \| \mathbf{r}_A \|_{A^{-1},  \, \Omega} \, \|\nabla e\|_{A, \, \Omega}  \mathrm{\:d}t & \leq
  \tfrac{1}{2}
    \int\limits_0^T \Big(
    {\alpha}_1(t) \| \mathbf{r}_A \|^2_{A^{-1}, \, \Omega} +
    \tfrac{1}{{\alpha}_1(t)} \, \|\nabla e\|^2_{A, \, \Omega}
    \Big) \mathrm{\:d}t,
    \label{eq:young-fenchel-for-decomposed-domain-2-ii-4} \\
\int\limits_0^T \left \| \tfrac{1}{\varrho} \, {\mathbf{r}_{f,\, \mu} } \right \|_{\Omega} \left \| \varrho \, e \right \|_{\Omega} \mathrm{\:d}t & \leq \tfrac{1}{2}
    \int\limits_0^T \Big(
    {\rho}_1(t) \, \left \| \tfrac{1}{\varrho} \,\mathbf{r}_{f,\, \mu} \right \|^2_{\Omega} +
    \tfrac{1}{{\rho}_1(t)} \,\left \| \varrho \, e \right \|^2_{\Omega}
    \Big) \mathrm{\:d}t,
    \label{eq:young-fenchel-for-decomposed-domain-2-ii-1}\\
%
    \int\limits_0^T R_{\,\mathcal{O}_{\mathrm{P}}, \{ \cdot \} }^{\footnotesize \frac{1}{\,2}} \|\varrho \, e \big \|_{\Omega} \mathrm{\:d}t & \leq \tfrac{1}{2}
  \int\limits_0^T \Big(
    {\rho}_2(t) R_{\,\mathcal{O}_{\mathrm{P}}, \{ \cdot \} } \,+\,
    \tfrac{1}{{\rho}_2(t)} \,\big \|\varrho \, e \big \|^2_{\Omega}
    \Big) \mathrm{\:d}t,
    \label{eq:young-fenchel-for-decomposed-domain-2-ii-2} \\
    %
    %
    %
    \int\limits_0^T R_{\,\mathcal{S}_{R}}^{\footnotesize \frac{1}{\,2}} \; \| \nabla e \|_{A, \, \Omega} \mathrm{\:d}t & \leq
    \tfrac{1}{2} \int\limits_0^T \Big(
    {\alpha}_3(t) R_{\,\mathcal{S}_{R}}(t) +
    \tfrac{1}{ {\alpha}_3(t)} \| \nabla e \|^2_{A, \, \Omega}
    \Big) \mathrm{\:d}t,
    \label{eq:young-fenchel-for-decomposed-domain-2-ii-6}
\end{alignat}
and
\begin{multline}
    \int\limits_0^T \big( R_{\,\mathcal{O}_{\mathrm{P}}, \, \| \cdot \|} + R_{\,\mathcal{O}_0}\big)^{\footnotesize \frac{1}{\,2}} \; \| \nabla e \|_{A, \Omega} \mathrm{\:d}t \\
    \leq \tfrac{1}{2} \int\limits_0^T \Big(
    {\alpha}_2(t) \big( R_{\,\mathcal{O}_{\mathrm{P}}, \, \| \cdot \| } + R_{\,\mathcal{O}_0}\big) \,+\,
    \tfrac{1}{{\alpha}_2(t)} \, \| \nabla e \|^2_{A, \, \Omega}
    \Big) \mathrm{\:d}t.
    \label{eq:young-fenchel-for-decomposed-domain-2-ii-3}
\end{multline}
%
By combining
(\ref{eq:young-fenchel-for-decomposed-domain-2-ii-4})--(\ref{eq:young-fenchel-for-decomposed-domain-2-ii-3}), we arrive at (\ref{eq:majorant-decomposed-i}).


\vskip 10pt

\noindent
(ii)
Existence of the pair $(v, y) \in V_0 \times Y_{\mathrm{div} \:}(Q_T)$ minimizing
the functional  \linebreak
$\overline{\mathrm M}^{\,2}_{\,\rm I, N} (v, y; \delta, \rho_1, \rho_2, \mu)$ can be proven
straightforwardly.
Indeed, let $v = u$ and $y = A \nabla u$. Since
$\mathrm{div} \: (A \nabla u) \in L^{2}(Q_T)$,
we see that $y \in Y_{\mathrm{div} \:}(Q_T)$. In this case (cf.
(\ref{eq:parabolic-problem-equation})--(\ref{eq:parabolic-problem-robin-boundary-condition})),
\begin{eqnarray*}
   & e(x, 0) &= \: (u - v)(x, 0) = u_0(x) - v(x, 0) = 0, \nonumber\\
   & \mathbf{r}_f(u, A \nabla u) &= \: f - u_t - \varrho^2 \, u + \mathrm{div} \: A \nabla u = 0, \nonumber\\
   & \mathbf{r}_A (u, A \nabla u) &= \: A \nabla u - A \nabla{u} = 0, \nonumber \\
   & \mathbf{r}_F (u, A \nabla u) &= \: F - \sigma^2 v - A \nabla u \cdot n = 0, \nonumber
\end{eqnarray*}
Thus, we see that $\overline{\mathrm M}_{\,\rm I, N} = 0$.  Since the majorant is nonnegative, the functions
$u$ and $A\nabla u$ minimize it.

Assume that ${\overline{\mathrm M}^{\,2}_{\,\mathrm{I}, \mathrm{N}}} = 0$. Then, the following relations hold:
\begin{alignat}{4}
    y = A \nabla v \;
		& \quad \mbox{a.a.} && \quad (x, t) \in Q_T,
    \label{eq:requirements-set-1}\\
    f - v_t - \varrho^2 \, v + \mathrm{div} \: y = 0 \;
		& \quad \mbox{a.a.} && \quad (x, t) \in \Omega_i \times (0, T), \quad
		                                        \Omega_i \subset \mathcal{O}_\Omega,
    \label{eq:requirements-set-2}\\
    v(\cdot, 0) = u_0                     \;
		& \quad \mbox{a.a.} && \quad  x \in \Omega,
    \label{eq:requirements-set-3}\\
    v = 0                       \;
		& \quad \mbox{a.a.} && \quad (x, t) \in S_D,
    \label{eq:requirements-set-4}\\
    y \cdot n + \sigma^2 v = F          \;
		& \quad \mbox{a.a.} && \quad (x, t) \in {\Gamma_R}_{j} \times (0, T), \quad
																						{\Gamma_R}_{j} \subset \mathcal{S}_{R}.
  \label{eq:requirements-set-5}
\end{alignat}
%
%
From (\ref{eq:requirements-set-2})--(\ref{eq:requirements-set-5}), it follows that
for any $\eta \in V_0$
\begin{equation*}
  \int\limits_0^T
    \sum\limits_{\Omega_i \subset \, \mathcal{O}_\Omega} \,
		\int\limits_{\Omega_i}
		\Big(( f - v_t - \varrho^2 \, v ) \, \eta  - y \cdot \nabla{\eta} \Big)\mathrm{\:d}x \,
    \,+\,
    \int\limits_0^T \sum\limits_{{\Gamma_R}_{j} \subset \, \mathcal{S}_{R}} \, 
             \int\limits_{{\Gamma_R}_{j}} F \eta {\mathrm{\:d}s\mathrm{d}t} = 0,
\end{equation*}
or, equally,
\begin{equation}
  \int\limits_{Q_T} \Big(( f - v_t - \varrho^2 \, v )\, \eta - y \cdot \nabla{\eta} ) \mathrm{\:d}x \, +\,
    \int\limits_{S_R} F \eta {\mathrm{\:d}s\mathrm{d}t} \,=\, 0,
    \quad \forall \eta \in V_0.
    \label{eq:generalized-statement-for-v-1}
\end{equation}
%
In view of (\ref{eq:requirements-set-1}), the identity
(\ref{eq:generalized-statement-for-v-1}) is
equivalent to (\ref{eq:generalized-statement}), whence it follows that $v = u$ and
$y = A \nabla u$. 

We
conclude that the exact lower bound of $\overline{\mathrm M}^{\,2}_{\,\rm I, N}$
is equal to zero and it is attained only on the pair $(v,y)$, which presents the exact solution
of
(\ref{eq:parabolic-problem-equation})--(\ref{eq:parabolic-problem-robin-boundary-condition})
and the respective flux.
\hfill $\Box$

\vskip 15pt


\subsection{Equivalence of $\overline{\mathrm M}_{\,\rm I, N}$ and the primal--dual error norm}
Now, we are aimed to  show that the majorant is equivalent to the error measure in terms
of a combined
(primal-dual) norm. This fact justifies the majorant as an adequate tool of error control.

Consider the solution of
(\ref{eq:parabolic-problem-equation})--(\ref{eq:parabolic-problem-robin-boundary-condition})
as a pair $(u, p) \in V_0 \times Y_{\mathrm{div} \:}(Q_T)$. In order to measure
the deviation of the approximation
$(v, y) \in V_0 \times Y_{\mathrm{div} \:}(Q_T)$ from $(u, p)$, we use the following
form of combined primal-dual norm
\begin{multline}
	\left \| [(u, p) - (v, y)] \right \|^2_{(\check{\nu}, \check{\theta}, \check{\zeta}, \check{\kappa}, \check{\chi}, \check{\vartheta}, \check{\varpi})} \\[-2pt]
	:=
	\check{\nu} \! \int\limits_0^T \! \|  \nabla{u - v} \|^2_{A, \, \Omega} \mathrm{\:d}t \,+\,
	\check{\kappa} \int\limits_0^T \! \| \varrho \, (u - v) \|^2_{\Omega} \mathrm{\:d}t \,+\,
	\check{\chi} \| (u - v) (\cdot, T) \|^2_{\Omega} \,\qquad \qquad \\[-10pt]
	+ \, \check{\theta} \! \int\limits_0^T \! \|  y - p \|^2_{A^{-1}, \, \Omega} \mathrm{\:d}t \,+\,
	\check{\zeta} \! \int\limits_0^T  \! \left \|  \mathrm{div} \: (p - y) - (u - v)_t \right\|^2_{\Omega} \mathrm{\:d}t \, \\[-10pt]
	+ \, \check{\vartheta} \int\limits_0^T \| \sigma\, (u - v) \|^2_{\Gamma_R} \mathrm{\:d}t \,+\,
	\check{\varpi}    \int\limits_0^T \| (p - y) \cdot n \|^2_{\Gamma_R} \mathrm{\:d}t.
	\label{eq:combined-norm}
\end{multline}
It is easy to see that the first three terms of (\ref{eq:combined-norm}) present an energy
norm of the error in the primal variable, the forth can be viewed as an error associated with
the flux. The fifth term is generated by both errors in primal and dual variables. The last two
terms are related to errors in boundary conditions.
%
For simplicity, further (\ref{eq:combined-norm}) is used as
$\left \| [(u, p) - (v, y)] \right \|^2$.

From Theorem \ref{th:theorem-minimum-of-majorant-I} (with $\alpha_1$, $\alpha_2$,
$\alpha_3 = {\rm const}$, $\mu = 0$, and exactly satisfied initial
condition $u_0 = v(\cdot, 0)$), the estimate can be written in the form
\begin{multline}
    (2 - \delta) \int\limits_0^T \|\, \nabla e \,\|^2_{A, \, \Omega} \,+\,
    \big(2 - \tfrac{1}{\gamma} \big) \int\limits_0^T \| \,\varrho\, e \,\|^2_{\Omega} \,+\,
    \| e (x, T) \|^2_{\Omega} \, \\[-10pt]
    + \, 2 \int\limits_0^T \|\, \sigma\, e \,\|^2_{\Gamma_R} \mathrm{\:d}t \,\leq\,
    {\overline{\mathrm M}^{\,2}_{\,\mathrm{I}, \mathrm{N}}} (v, y) \, :=\,
    \gamma \int\limits_0^T R_{\mathcal{O}_{\mathrm{P}}, \, \{ \cdot \}} \mathrm{\:d}t \, \qquad \qquad \qquad \qquad \\[-10pt]
    + \,
    \alpha_1 \int\limits_0^T \| \, \mathbf{r}_A \, \|^2_{A^{-1}, \Omega} \mathrm{\:d}t \,+\,
    \alpha_2 \int\limits_0^T \Big( \! R_{\mathcal{O}_{\mathrm{P}}, \, \| \cdot \|} \,+\,
                               R_{\mathcal{O}_{0}} \Big) \mathrm{\:d}t \, + \,
    \alpha_3 \int\limits_0^T \! R_{\mathcal{S}_{R}} \mathrm{\:d}t.
    \label{eq:majorant-simplified-combined-norm}
\end{multline}
%
Set
\begin{alignat}{2}
C_{\Omega_l \mathrm{P}} &
:= \max\limits_{ \Omega_l \subset \,\mathcal{O}_\mathrm{P}}
   \left\{ \, \tfrac{|\Omega_l|}{\!\mathrm{P}^2} \right\}, \qquad \qquad \;
\overline{C}_{\mathrm{P\Omega}}
:= \max\limits_{ \Omega_i \subset \,\mathcal{O}_\Omega } \Big\{ \, {C_{\mathrm{P\Omega}}}_i \Big\}, \\
\overline{C}_{\mathrm{\Gamma\Omega}} &
:= \max\limits_{ {\Gamma_R}_{j} \subset \, \mathcal{S}_{R} } \Big\{ \, {C_{\mathrm{\Gamma\Omega}}}_{j} \, \Big\},
 \qquad \quad
C_{\gamma \alpha_2} := \max \{ \gamma, \alpha_2\},
\end{alignat}
then
\begin{multline}
	\overline{\mathrm M}^{\,2}_{\,\rm I, N} \,
	\leq\,  %
	\| e (\cdot, T) \|^2_{\Omega} \,
	+ \, \alpha_1 \int\limits_0^T \| \, y - A \nabla{v} \, \|^2_{A^{-1}, \Omega}\mathrm{\:d}t \, \\[-5pt]
	+ \tfrac{C_{\gamma \alpha_2}}{\!\underline{\lambda}_A}
		\max \Big\{ \overline{C}_{\mathrm{P\Omega}}^2, C_{\Omega_l \mathrm{P}} \Big \}
		\int\limits_0^T  \| \, f - v_t + \mathrm{div} \: y - \varrho^2 v \, \|^2_{\Omega} \mathrm{\:d}t \,\\[-10pt]
	+ \, \alpha_3 \, \tfrac{\overline{C}_{\mathrm{\Gamma\Omega}}^2}{\!\underline{\lambda}_A}
		\int\limits_0^T \| \, F - \sigma^2 v - y \cdot n \, \|^2_{\Gamma_R} \mathrm{\:d}t.
	\label{eq:majorant-estimated-with-max-constants}
\end{multline}
%
For further simplification, let
\begin{equation}
C_{\rm max} :=
\tfrac{C_{\gamma \alpha_2}}{\!\underline{\lambda}_A}
\max \Big\{ \overline{C}_{\mathrm{P\Omega}}^2, C_{\Omega_l \mathrm{P}} \Big \}, \qquad
C_{\rm \alpha_3 \Gamma} :=
\alpha_3 \, \tfrac{\overline{C}_{\mathrm{\Gamma\Omega}}^2}{\!\underline{\lambda}_A}.
\label{eq:c-max-c}
\end{equation}
By means of (\ref{eq:parabolic-problem-equation}) and
(\ref{eq:parabolic-problem-robin-boundary-condition}), the right-hand side of
(\ref{eq:majorant-estimated-with-max-constants}) can be decomposed as follows:
\begin{multline}
    {\overline{\mathrm M}^{\,2}_{\,\mathrm{I}, \mathrm{N}}} \,\leq\,  %
    \| e (\cdot, T) \|^2_{\Omega} \,+\,
    \alpha_1 \left( \int\limits_0^T \|  y - p \|^2_{A^{-1}, \, \Omega} \mathrm{\:d}t \,+\,
                    \int\limits_0^T \|  \nabla (u - v)\|^2_{A, \, \Omega} \mathrm{\:d}t \right) \,
    \qquad \qquad \qquad \qquad \qquad \qquad \qquad \quad\\
    +\, C_{\rm max}
    \Big( \int\limits_0^T \| \, \mathrm{div} \: (y - p) + (u - v)_t \, \|^2_{\Omega} \mathrm{\:d}t \,+\,
    \int\limits_0^T \Big \| \varrho^2 (u - v) \Big \|^2_{\Omega} \mathrm{\:d}t \Big) \, \\
    +\, C_{\rm \alpha_3 \Gamma}
    \left( \int\limits_0^T \Big \| \sigma^2 (u - v) \Big \|^2_{\Gamma_R} \mathrm{\:d}t \,+\,
    \int\limits_0^T \| (p - y) \cdot n \|^2_{\Gamma_R} \mathrm{\:d}t \right). \nonumber
\end{multline}
We recall that $\varrho$ and $\sigma$ are uniformly bounded by the constants
$C_{\mathrm{\varrho}}$ and $C_{\mathrm{\sigma}}$, respectively,  and
we estimate the right-hand side of the latter inequality as
\begin{multline}
	{\overline{\mathrm M}^{\,2}_{\,\mathrm{I}, \mathrm{N}}} \leq  %
	\| e (\cdot, T) \|^2_{\Omega} \,+\,
	\alpha_1
	\left( \int\limits_0^T \|  y - p \|^2_{A^{-1}, \, \Omega} \mathrm{\:d}t \, + \int\limits_0^T \|  \nabla (u - v)\|^2_{A, \, \Omega} \mathrm{\:d}t \right) \, \\
  + C_{\rm max} \left( \int\limits_0^T \| \mathrm{div} \: (y - p) + (u - v)_t \|^2_{\Omega} \mathrm{\:d}t \,+\,
	 C_{\mathrm{\varrho}}^2 \int\limits_0^T \| \varrho\, (u - v) \|^2_{\Omega_{\mathrm{P}}} \mathrm{\:d}t \right) \,
	\qquad\\
	\qquad \qquad \qquad
	+\, C_{\rm \alpha_3 \Gamma} \left(
	C_{\mathrm{\sigma}}^2 \int\limits_0^T \| \sigma (u - v) \|^2_{\Gamma_R} \mathrm{\:d}t + \int\limits_0^T \| (p - y) \cdot n \|^2_{\Gamma_R} \mathrm{\:d}t \right) \\
	\qquad \qquad \qquad \qquad
	=: \big \| [(u, p) - (v, y)] \big \|^2_{(\check{\nu}, \check{\theta}, \check{\zeta}, \check{\kappa}, \check{\chi}, \check{\vartheta}, \check{\varpi})}.
	\label{eq:majorant-estimated-with-lambda-sigma-constants}
\end{multline}
Here, on the right-hand side we have the error measured in terms of
(\ref{eq:combined-norm}) with positive weights
\begin{equation}
	\check{\nu} = \check{\theta} = \alpha_1, \quad
	\check{\zeta} = C_{\rm max}, \quad
	\check{\kappa} = C_{\mathrm{\varrho}}^2 \,C_{\rm max}, \quad
	\check{\chi} = 1, \quad
	\check{\vartheta} = C_{\mathrm{\sigma}}^2 \, C_{\rm \alpha_3 \Gamma}, \quad
	\check{\varpi} = C_{\rm \alpha_3 \Gamma}.
\end{equation}
%
%
Next, we combine four terms related to the energy error norm of the primal variable
on the right-hand side of
(\ref{eq:majorant-estimated-with-lambda-sigma-constants}) and estimate it by
using (\ref{eq:majorant-simplified-combined-norm}). The rest of the terms related to
dual component can be estimated by the technique used above. Therefore, we obtain
\begin{multline}
    \left \| [(u, p) - (v, y)] \right \|^2 \leq
    C_{\rm ER} \, {\overline{\mathrm M}^{\,2}_{\,\mathrm{I}, \mathrm{N}}} \,+\,
    \alpha_1 \left( \int\limits_0^T \|  y - A \nabla v \|^2_{A^{-1}, \Omega} \mathrm{\:d}t \,+\, \int\limits_0^T \|  \nabla (u - v)\|^2_{A, \Omega} \mathrm{\:d}t  \right) \, \\
    + \,C_{\rm max} \bigg(
    \int\limits_0^T \| f + \mathrm{div} \: y - v_t  - \varrho^2 v\|^2_{\Omega}\mathrm{\:d}t \,+\,
    C_{\mathrm{\varrho}}^2 \int\limits_0^T \| \varrho\, (u - v) \|^2_{\Omega}\mathrm{\:d}t \bigg) \, \\
    + \, C_{\rm \alpha_3 \Gamma}
		\left( \int\limits_0^T \| F - \sigma^2 v - y \cdot n\|^2_{\Gamma_R} \mathrm{\:d}t
    \,+\,
    C_{\mathrm{\sigma}}^2 \int\limits_0^T \|\sigma (v - u) \|_{\Gamma_R} \mathrm{\:d}t \right),
    \label{eq:estimate-primal-error-norm}
\end{multline}
where
\begin{equation}
	C_{\rm ER} = \max \left\{
	\tfrac{\alpha_1}{(2 - \delta)}, \quad
	\tfrac{\gamma}{(2 \gamma - 1)} \,C_{\mathrm{\varrho}}^2 \, C_{\rm max}, \quad
	1, \quad
	\tfrac12 C_{\mathrm{\sigma}}^2 \,C_{\rm \alpha_3 \Gamma}
	\right\}.
\end{equation}
By using constants
\begin{equation*}
\underline{C}_{\mathrm{P\Omega}} := \min\limits_{ \Omega_i \subset \, \mathcal{O}_\Omega} \left\{ \, {C_{\mathrm{P\Omega}}}_i \, \right\}, \quad \mbox{and} \quad
\widetilde{C}_{\mathrm{\Gamma\Omega}} := \tfrac{\overline{C}_{\mathrm{\Gamma\Omega}}^2}{\!\underline{C}_{\mathrm{\Gamma\Omega}}^2}, \quad \mbox{where}
\quad
\underline{C}_{\mathrm{\Gamma\Omega}} := \min\limits_{ {\Gamma_R}_{j} \subset \, \mathcal{S}_{R}} \left\{ \, {C_{\mathrm{\Gamma\Omega}}}_{j} \, \right\},
\end{equation*}
we rewrite the right-hand side of (\ref{eq:estimate-primal-error-norm})
and obtain the following result:
%
%
\begin{multline}
\left \| [(u, p) - (v, y)] \right \|^2 \leq
\alpha_1 \int\limits_0^T  \| \, \nabla (u - v) \,\|^2_{A, \, \Omega} \mathrm{\:d}t \, +\,
C_{\mathrm{\varrho}}^2 \, C_{\rm max}
\int\limits_0^T \| \, \varrho\, (u - v) \, \|^2_{\Omega} \mathrm{\:d}t \,
+\,\| e (\cdot, T) \|^2_{\Omega} \, \qquad \qquad \qquad \\
+\,
C_{\mathrm{\sigma}}^2 \, C_{\rm \alpha_3 \Gamma}
\int\limits_0^T \| \, \sigma (v - u) \, \|^2_{\Gamma_R} \mathrm{\:d}t \,+\,
C_{\rm ER} \, {\overline{\mathrm M}^{\,2}_{\,\mathrm{I}, \mathrm{N}}} \,+\,
\alpha_1 \int\limits_0^T \| \, y - A \nabla{v} \,\|^2_{A^{-1}, \Omega} \mathrm{\:d}t \,\\
+\, \, \tfrac{ C_{\rm max} }{\underline{C}_{\mathrm{P\Omega}}^2}
\left( \int\limits_0^T R_{\mathcal{O}_{\mathrm{P}}, \, \| \cdot \|} \mathrm{\:d}t +
\int\limits_0^T R_{\mathcal{O}_0} \mathrm{\:d}t \right) \, + \,
\alpha_3 \, \widetilde{C}_{\mathrm{\Gamma\Omega}} \int\limits_0^T R_{\mathcal{S}_{R}} \mathrm{\:d}t\,.
\label{eq:estimate-primal-error-norm-2}
\end{multline}
Finally, the terms related to the error norm of primal component on the right-hand
side of (\ref{eq:estimate-primal-error-norm-2})
can be estimated by the majorant (\ref{eq:majorant-simplified-combined-norm}):
\begin{multline*}
\big \| [(u, p) - (v, y)] \big \|^2 \leq
\int\limits_0^T \bigg (
\alpha_1 \left(2 C_{\rm ER} + 1 \right) \!\|  y - A \nabla{v} \|^2_{A^{-1}, \Omega} \,+\,
2 \, C_{\rm ER} \, \gamma R_{\mathcal{O}_{\mathrm{P}}, \, \{ \cdot \}} \,+\, \\
\qquad \qquad \qquad \qquad
\alpha_2 \left( 2\, C_{\rm ER} + \tfrac{C_{\rm max} \,}{\alpha_2 \underline{C}_{\mathrm{P\Omega}}^2} \right)
\left (R_{\mathcal{O}_{\mathrm{P}}, \, \| \cdot \|} \,+\, R_{\mathcal{O}_0} \right)
%
%
\, \\ + \,
\alpha_3 \, \left(2\, C_{\rm ER} + \widetilde{C}_{\mathrm{\Gamma\Omega}} \right) R_{\mathcal{S}_{R}}
\bigg ) \mathrm{\:d}t \leq C_{\rm MAJ} \, {\overline{\mathrm M}^{\,2}_{\,\mathrm{I}, \mathrm{N}}},
\end{multline*}
where
\begin{equation}
C_{\rm MAJ} = \max \left\{
2\, C_{\rm ER} + 1, \quad
2\, C_{\rm ER}, \quad
2\, C_{\rm ER} +
\tfrac{C_{\rm max}\,}{\!\alpha_2 \underline{C}_{\mathrm{P\Omega}}^2}, \quad
2\, C_{\rm ER} + \widetilde{C}_{\mathrm{\Gamma\Omega}}
\right\}.
\end{equation}
Therefore, we obtain the double inequality
\begin{equation}
{\overline{\mathrm M}^{\,2}_{\,\rm I, N}} \! \leq \!
\big \| [(u, p) - (v, y)] \big \|^2_{(\check{\nu}, \check{\theta}, \check{\zeta}, \check{\kappa}, \check{\chi}, \check{\vartheta})} \leq C_{\rm MAJ} \, {\overline{\mathrm M}^{\,2}_{\,\rm I, N}},
\label{eq:double-inequality-maj-I}
\end{equation}
which shows that the majorant introduced in Theorem
\ref{th:theorem-majorant-for-decomposed-domain-1} is equivalent to a certain form of combined
(primal-dual) error norm. In other words, 
${\overline{\mathrm M}^{\,2}_{\,\rm I, N}}(v, y; \delta, \rho_1, \rho_2, \mu)$ 
(which contains only known functions and
parameters) adequately reflects the distance between \linebreak
$(v, y) \in V_0 \times Y_{\mathrm{div} \:}(Q_T)$ and the exact solution $(u, p)$.
In particular, this means that if $(u_h, p_h)$ is the sequence of approximations
computed on a certain set of meshes $\mathcal{F}_h$, which converges to $(u, p)$
with the rate $h^\alpha$, then the values of the majorant
tend to zero with the same rate.


\subsection{The second estimate}
Now, we deduce another  estimate, which is in general sharper than
(\ref{eq:majorant-decomposed-i}), but contains an additional free function
$w\in V_0$. The corresponding residuals of (\ref{eq:dual-part}),
(\ref{eq:parabolic-problem-equation}), and
(\ref{eq:parabolic-problem-robin-boundary-condition}) are presented as
\begin{alignat}{2}
\mathbf{r}_f (v, y, w) & := f - {(v + w)}_t - \varrho^2 \, (v - w) + \mathrm{div} \: y, \label{eq:r-1} \\
\mathbf{r}_{f,\, \mu} (v, y, w) & := \mu \,\mathbf{r}_f (v, y, w), \label{eq:r-1-mu} \\
\mathbf{r}_{f, 1 - \mu} (v, y, w) & := (1 - \mu) \, \mathbf{r}_f (v, y, w), \label{eq:r-1-1-mu} \\
\mathbf{r}_A (v, y, w) & := y - A \nabla{(v - w)}, \label{eq:r-2} \\
\mathbf{r}_F (v, y, w) & := F - \sigma^2 (v - w) - y \cdot n, \label{eq:r-3}
\end{alignat}
respectively.
On collections $\mathcal{O}_0$ and $\mathcal{S}_{R}$, we impose the mean conditions
similar to (\ref{eq:mean-condition-on-omega-0}) and (\ref{eq:mean-condition-on-boundary}),
namely,
\begin{equation}
    { \Big \{ \mathbf{r}_{f, 1 - \mu} (v, y, w) \Big\} }_{{\Omega}_k \subset\, \mathcal{O}_0} = 0,
		\quad
    \mbox{for a.a.} \quad t \in [0, T],
    \label{eq:r-1-mean-condition-on-omega-0}
\end{equation}
and
\begin{equation}
    { \Big \{ \mathbf{r}_F (v, y, w) \Big\} }_{{\Gamma_R}_{j} \subset \mathcal{S}_{R}} = 0, \quad
    \mbox{for a.a.}\quad t \in [0, T].
    \label{eq:r-3-mean-condition-on-omega-0}
\end{equation}
Correspondingly, the complexes
$R_{\, \mathcal{O}_{\mathrm{P}}, \{ \mathbf{r}_{f, 1 - \mu} \}}(t)$,
$R_{\, \mathcal{O}_{\mathrm{P}}, \| \mathbf{r}_{f, 1 - \mu} \|}(t)$,
$R_{\,\mathcal{O}_0}(t)$, $R_{\mathcal{S}_{R}}(t)$ are defined analogously
(\ref{eq:r-op-mean})--(\ref{eq:r-sr}) and depend on
residuals (\ref{eq:r-1})--(\ref{eq:r-3}), which are based on free functions $v$, $y$, 
and $w$.
%
%
\begin{theorem}
\label{th:theorem-decomposed-domain-majorant-II-1}
(i)
Assume that conditions (\ref{eq:r-1-mean-condition-on-omega-0}) and
(\ref{eq:r-3-mean-condition-on-omega-0}) are satisfied.
Then, for any \linebreak
$v, w \in V_0$ and $y \in Y_{\mathrm{div} \:}(Q_T)$,
${\delta} \in (0, 2]$, ${\epsilon} \geq 1$, ${\rho}_1(t) \geq 1$,
${\rho}_2(t) \geq 1$,
the error has the following estimate:
\begin{multline}
    {[e]}^{\,2}_{({\nu},\, {\theta},\, {\zeta}, \, 2)} \leq
    \overline{\mathrm M}^{\,2}_{\,\rm II, N} (v, y, w; \, \delta, \epsilon, \rho_1, \rho_2, \mu) :=
    {\epsilon} \| w(x, T)\|^2_{\Omega} \, + \,
    2L(v, w) \, + \, l(v, w) \, \\
    + \quad \int\limits_0^T \! \! \bigg( \!
    {\rho}_1(t) \Big\|\tfrac{1}{\varrho} \, {\mathbf{r}_{f,\, \mu} (v, y ,w)} \Big\|^2_{\Omega_{\mathrm{P}}} \! \! +
    {\rho}_2(t) R_{\, \mathcal{O}_{\mathrm{P}}, \{ \cdot \}}(t) \, + \,
    {\alpha}_1(t) \| \mathbf{r}_A (v, y ,w)\|^2_{A^{-1}, \, \Omega} \,\\[-8pt]
    + {\alpha}_2(t) \left( R_{\, \mathcal{O}_{\mathrm{P}}, \| \cdot \|}(t) \, + \,
     R_{\,\mathcal{O}_0}(t) \right) \,+\,
    {\alpha}_3(t) R_{\,\mathcal{S}_{R}}(t)
    \bigg) \! \! \mathrm{\:d}t,
    \label{eq:majorant-2-decomposed-domain-ii}
\end{multline}
\begin{equation}
L(v, w) :=
\int\limits_{Q_T} \Big( v_t \,w + A \nabla v \cdot \nabla w + \varrho^2 \, v \,w  - f w \Big) {\mathrm{\:d}x\mathrm{d}t} \,-\,
\int\limits_{S_R} (F - \sigma^2 v) \,w {\mathrm{\:d}s\mathrm{d}t},
\label{eq:l-function}	
\end{equation}
\begin{equation}
l(v, w) := \int\limits_\Omega |v(x, 0) - u_0(x)|^2 - 2 w(x, 0) \big( u_0(x) - v(0, x)\big) \mathrm{\:d}x,
\label{eq:l-small-function}	
\end{equation}
%
%
%
${\nu} = 2 - {\delta}$,
${\theta}(x, t) =
\varrho(x) \left(2 - \tfrac{1}{{\rho_1(t)}} - \tfrac{1}{{\rho_2(t)}}
           \right)^{\footnotesize \frac{1}{2}}\,$,
$\zeta = 1 - \tfrac{1}{\epsilon}$ are positive parameters 
,
$\mu(x,t)\in [0, 1]$ is real-valued function,
${\alpha}_1(t)$, ${\alpha}_2(t)$, ${\alpha}_3(t)$ are positive functions satisfying
(\ref{eq:alpha}).
\vskip 10pt
\noindent
(ii)
For any $\delta \in (0, 2]$, ${\rho}_1(t) \geq 1$,
${\rho}_2 \geq 1$, ${\epsilon} \geq 1$, and ${\mu} \in [0, 1]$, the lower bound of the
variation problem generated by the majorant
\begin{equation}
    \inf\limits_{
    \begin{array}{c}
    v, w \in V_0\\
    y \in Y_{\mathrm{div} \:}(Q_T)
    \end{array}
    } {\overline{\mathrm M}^{\,2}_{\,\mathrm{II}, \mathrm{N}}}(v, y, w; \, \delta, \epsilon, \rho_1, \rho_2, \mu)
    \label{eq:inf-maj-II}
\end{equation}
is zero, and it is attained if and only if $v = u$, $y = A \nabla u$, and $w = 0$.
\end{theorem}

\noindent
{\bf Proof:}
(i)
We rewrite the right-hand side of (\ref{eq:energy-balance-equation}) by inserting
functions $w \in V_0$ and $y \in Y_{\mathrm{div} \:}(Q_T)$, which implies the following
relation
\begin{multline}
    \int\limits_0^T \| \nabla e \|^2_{A, \, \Omega} \mathrm{\:d}t \,+\,
    \int\limits_0^T \| \varrho \, e \|^2_{\Omega}\mathrm{\:d}t \,+\,
    \tfrac{1}{2} \| e(\cdot, T) \|^2_{\Omega} \,+\,
    \int\limits_0^T \| \sigma\, e \|^2_{\Gamma_R} \mathrm{\:d}t \,\\
    = \int\limits_{\Omega} e(x, T)w(x, T) \mathrm{\:d}x +
    \int\limits_{\Omega} \Big(\tfrac{1}{2} \, e^2(x, 0) - e(x, 0)w(x, 0) \Big)\mathrm{\:d}x\, \\
    + \mathscr{I}_1 + \mathscr{I}_2 + \mathscr{I}_3 + L(v, w),
    \label{eq:decomposition-for-majorant-II}
\end{multline}
where $\mathscr{I}_f$, $\mathscr{I}_A$, $\mathscr{I}_F$ are quite analogous to (\ref{eq:Ir-Id-Ib-terms})
and depend on residuals (\ref{eq:r-1})--(\ref{eq:r-3}).
%
%
%
Following the steps of the proof of Theorem \ref{th:theorem-majorant-for-decomposed-domain-1},
the term $\mathscr{I}_f$ can be represented as
\begin{equation}
    \mathscr{I}_f =
    \mathscr{I}_{f, \,\mu} \, + \, {\mathscr{I}}^{\mathcal{O}_{\mathrm{P}}}_{f, \,1 - \mu} \, + \,
    {\mathscr{I}}^{\mathcal{O}_0}_{f, \,1 - \mu},
    \label{eq:representation-of-if-with-mu}
\end{equation}
where each of the summands is estimated as follows
\begin{alignat}{2}
\mathscr{I}_{f, \mu}     &
\leq \int\limits_0^T \left \| \tfrac{1}{\varrho} \,{\mathbf{r}_{f,\, \mu}} \right \|_{\Omega}
              \left \| \varrho \, e \right \|_{\Omega} \mathrm{\:d}t,
\label{eq:i-mu-majorant-on-decomposed-omega}\\
{\mathscr{I}}^{\mathcal{O}_{\mathrm{P}}}_{f, 1 - \mu} &
\leq
\int\limits_0^T R^{\tfrac{1}{\,2}}_{\,\mathcal{O}_{\mathrm{P}}, \| \cdot \|} \;
\| \nabla e \|_{A, \, {\Omega}_{\mathrm{P}}} \mathrm{\:d}t \, +
\int\limits_0^T R^{\tfrac{1}{\,2}} _{\,\mathcal{O}_{\mathrm{P}}, \{ \cdot \}}
\big \| \varrho \, e \big \|_{\Omega} \mathrm{\:d}t,
\label{eq:i-1-mu-majorant-on-decomposed-omega}\\
{\mathscr{I}}^{\mathcal{O}_0}_{f, \,1 - \mu} & \leq
\int\limits_0^T R_{\, \mathcal{O}_0}^{\tfrac{1}{\,2}} \;
\| \nabla e \|_{A, \Omega_0} \mathrm{\:d}t.
\label{eq:r-oo-majorant-on-decomposed-omega}
\end{alignat}
%
%
The term $\mathscr{I}_A$ is estimated by using the H\"older inequality, and $\mathscr{I}_F$ 
is bounded
analogously to (\ref{eq:estimation-i-b}).
%
%
Following steps of the proof of Theorem \ref{th:theorem-majorant-for-decomposed-domain-1}, we use inequality
\begin{alignat}{2}
    \int\limits_\Omega e(x, T) \,w(x, T)  {\mathrm{\:d}x\mathrm{d}t} & \leq
    \tfrac{1}{2}
		\Big(
		\tfrac{1}{\epsilon} \, \| e(\cdot, T)\|^2_{\Omega} +
    \epsilon \| w(\cdot, T)\|^2_{\Omega}
		\Big), \label{eq:ineq-young-fenchel-1}
    %
		\end{alignat}
to estimate the term related to $t=T$, and the Young--Fenchel inequalities to estimate
(\ref{eq:i-mu-majorant-on-decomposed-omega})--(\ref{eq:r-oo-majorant-on-decomposed-omega})
, $\mathscr{I}_A$, and $\mathscr{I}_F$.
By combining these estimate all together,
we obtain the required estimate (\ref{eq:majorant-2-decomposed-domain-ii}).
\vskip 15pt

\noindent
(ii)
The proof is similar to the proof of (ii) in
Theorem \ref{th:theorem-majorant-for-decomposed-domain-1}.
\hfill $\Box$

\vskip 15pt
\subsection{Equivalence of $\overline{\mathrm M}_{\,\rm II, N}$ and 
$[ e ]_{(\nu, \theta, \zeta, \chi)}$}
Finally, we prove that ${\overline{\mathrm M}^{\,2}_{\,\mathrm{II}, \mathrm{N}}}$ is equivalent to the error measure
(\ref{eq:energy-norm-for-reaction-diff-evolutionary-problem}). For this purpose, we
estimate (\ref{eq:majorant-2-decomposed-domain-ii}) from above and show that this
upper bound is equivalent to the error norm. Henceforth, we assume that $\mu = 0$ (this
is done for the sake of simplicity only), $y = A \nabla u \in Y_{\mathrm{div} \:}(Q_T)$, and
$w = e$, then
\begin{equation*}
    \mathbf{r}_f (v, A \nabla u, e) = 2 \varrho^2\, e, \quad
    \mathbf{r}_A (v, A \nabla u, e) = 2 A \nabla e, \quad
    \mathbf{r}_F (v, A \nabla u, e) = 2 \sigma^2 e.
    \label{eq:r-1-r-2-r-3-changed}
\end{equation*}
%
%
The functional (\ref{eq:l-function}) can be represented as follows:
\begin{alignat}{2}
	L(v, e) = &
	\int\limits_{Q_T}
	\Big(  v_t \,e + A \nabla v \cdot \nabla e + \varrho^2\, v \, e - f e \Big) {\mathrm{\:d}x\mathrm{d}t}
	- \int\limits_{S_R} ( F - \sigma^2\, v) \,e {\mathrm{\:d}s\mathrm{d}t} \nonumber\\
	%
	= & \int\limits_{Q_T}
	\Big( u_t e + A \nabla u \cdot \nabla e + \varrho^2\, u e  - f e \Big) {\mathrm{\:d}x\mathrm{d}t}
	+ \int\limits_{S_R} ( F - \sigma^2\, u) \,e {\mathrm{\:d}s\mathrm{d}t} \nonumber\\
	& \qquad \qquad \qquad
	- \int\limits_{Q_T} \big(  A \nabla e \cdot \nabla e + e_t e + \varrho^2\, e^2 \big) {\mathrm{\:d}x\mathrm{d}t}
	- \int\limits_{S_R} \sigma^2 \,e^2 {\mathrm{\:d}s\mathrm{d}t}.
	\label{eq:l-function-1}
\end{alignat}
In view of (\ref{eq:generalized-statement}), the first two terms in the right-hand side
of (\ref{eq:l-function-1}) vanishes, and we find that
\begin{equation}
    L(v, e) =
		- \int\limits_{Q_T} \big(  A \nabla e \cdot \nabla e + e_t e + \varrho^2\, e^2 \big) {\mathrm{\:d}x\mathrm{d}t}
		- \int\limits_{S_R} \sigma^2 \,e^2 {\mathrm{\:d}s\mathrm{d}t}.
\end{equation}
Next,
\begin{equation}
    l(v, e) =
    \int\limits_\Omega \left( |v(x, 0) - u_0(x)|^2 - 2 e(x, 0) \big( u_0(x) - v(0, x)\big) \right)\mathrm{\:d}x =
    - \| e (x, 0)\|_{\Omega}^2.
    \label{eq:l-small-function-1}
\end{equation}
%
By means of differentiation by part and (\ref{eq:l-small-function-1}), we obtain the
estimate
\begin{multline}
    {\overline{\mathrm M}^{\,2}_{\,\rm II, N}} \leq
    \left ( 4 \,\alpha_2 - 2 \right) \int\limits_0^T \| \nabla e \|^2_{A, \, \Omega} \mathrm{\:d}t +
    4 \,\alpha_1 \int\limits_0^T \sum\limits_{\Omega_l \subset \mathcal{O}_{\mathrm{P}}} \tfrac{{C_{\mathrm{P\Omega}}}_l^2}{\!\underline{\lambda}_A}  \left\| \, \varrho\, e \,\right\|^2_{\Omega_l} \mathrm{\:d}t
    - 2 \int\limits_0^T \| \varrho\, e\|^2_{\Omega} \mathrm{\:d}t \\
    + 4 \,\alpha_3 \int\limits_0^T \sum\limits_{{\Gamma_R}_{j} \subset \, \mathcal{S}_{R}}
    \tfrac{{C_{\mathrm{\Gamma\Omega}}}_{j}^2}{\!\underline{\lambda}_A}  \left\| \sigma\, e \,\right\|^2_{{\Gamma_R}_{j_m}} \mathrm{\:d}t 
    - 2 \int\limits_0^T \| \sigma\, e\|^2_{\Gamma_R} \mathrm{\:d}t \\
		\qquad \qquad \qquad \qquad \qquad \qquad \qquad 
		+ \epsilon \| e(\cdot, T)\|^2_{\Omega} - 2 \int\limits_{Q_T} e_t e {\mathrm{\:d}x\mathrm{d}t}
		- \| e (x, 0)\|_{\Omega}^2\\
     \leq 2 \left ( 2 \,\alpha_2 - 1 \right) \int\limits_0^T \| \nabla e \|^2_{A, \, \Omega} \mathrm{\:d}t +
    2 \left ( 2 \,\alpha_1 \tfrac{\overline{C}_{\mathrm{P\Omega}}^2}{\!\underline{\lambda}_A} - 1 \right) \int\limits_0^T \| \varrho\, e \|^2_{\Omega}  \qquad \qquad \qquad \qquad \\
    + 2 \left ( 2 \,\alpha_3 \tfrac{\overline{C}_{\mathrm{\Gamma\Omega}}^2}{\!\underline{\lambda}_A} - 1 \right) \int\limits_0^T \left \| \sigma\, e \right\|^2_{\Gamma_R} +
    (\epsilon - 1) \| e(\cdot, T)\|^2_{\Omega},
\end{multline}
where
$\overline{C}_{\mathrm{P\Omega}} = \max\limits_{ \Omega_l \subset \mathcal{O}_{\mathrm{P}} } \left\{ \, {C_{\mathrm{P\Omega}}}_l \right\}$,
$\overline{C}_{\mathrm{\Gamma\Omega}} = \max\limits_{ {\Gamma_R}_j \subset \mathcal{}_R} \{ {C_{\mathrm{\Gamma\Omega}}}_{j} \}$.
Therefore, for any $v \in V_0$ we arrive at double inequality
\begin{equation}
{[e]}^{\,2}_{(\nu,\, \theta,\, \zeta, \, \chi)} \leq
{\overline{\mathrm M}^{\,2}_{\,\rm II, N}} \leq
{[e]}^{\,2}_{({\nu}',\, {\theta}',\, {\zeta}', \,  {\chi}')} \leq
\mathcal{K} {[e]}^{\,2}_{(\nu,\, \theta,\, \zeta, \, \chi)},
\label{eq:double-inequality-maj-II}
\end{equation}
with parameters
\begin{eqnarray*}
&&\nu' =  2 \left ( 2 \,\alpha_2 - 1 \right), \quad
\theta' =
\varrho \bigg( 2 \left( 2 \,\alpha_1 
\tfrac{\overline{C}_{\mathrm{P\Omega}}^2}
{\!\underline{\lambda}_A} - 1 \right) \bigg)^{\footnotesize \frac{1}{\,2}}, \\\
&&\zeta' = \epsilon - 1, \quad
\chi' = 2 \left ( 2 \,\alpha_3 \tfrac{\overline{C}_{\mathrm{\Gamma\Omega}}^2}{\!\underline{\lambda}_A} - 1 \right),\\
&&\nu = 2 - {\delta}, \quad
\theta = \varrho \left( 2 - \tfrac{1}{\,{\gamma}} \right)^{\footnotesize \frac{1}{\,2}}, \quad
\zeta = 1 - \tfrac{1}{\, {\epsilon}}, \quad
\chi = 2,
\end{eqnarray*}
and
\begin{equation*}
\mathcal{K} = \max
\Bigg \{
\tfrac{2\, (2 \,\alpha_2 - 1)}{2 - \delta}, \quad
2\, \Bigg(\tfrac{2 \,\alpha_1 \tfrac{\overline{C}_{\mathrm{P\Omega}}^2}{\!\underline{\lambda}_A} - 1}{2 - \tfrac{1}{\,{\gamma}}}
\Bigg)^{\footnotesize \frac{1}{\,2}}, \quad
\epsilon, \quad
2 \,\alpha_3 \tfrac{\overline{C}_{\mathrm{\Gamma\Omega}}^2}{\!\underline{\lambda}_A} - 1
\Bigg \}.
\end{equation*}
The  relation (\ref{eq:double-inequality-maj-II}) shows that
${\overline{\mathrm M}^{\,2}_{\,\rm II, N}}$ is
equivalent to the error measure (\ref{eq:energy-norm-for-reaction-diff-evolutionary-problem}).
Therefore, we obtain fully 
error majorants (presented in Theorems \ref{th:theorem-majorant-for-decomposed-domain-1} and
\ref{th:theorem-decomposed-domain-majorant-II-1}), which
generate fully computable and realistic estimates of the distance to exact solution.

\noindent
{\bf Acknowledgment}. In part, the work is supported by RFBR grant 11-01-00531-a.

\bibliographystyle{plain}
\def\cprime{$'$}

\end{document}